\documentclass[leqno]{article}

\usepackage{amssymb}

\newtheorem{theorem}{Theorem}[section]  
\newtheorem{thm}[theorem]{Theorem}  
\newtheorem{lm}[theorem]{Lemma}  
\newtheorem{prop}[theorem]{Proposition} 
\newtheorem{cor}[theorem]{Corollary}  

\newtheorem{myremark}[theorem]{Remark}
\newenvironment{rem}{\begin{myremark}\rm}{\end{myremark}}

\newtheorem{myexample}[theorem]{Example}

\newtheorem{mydefinition}[theorem]{Definition}

\def\beq{\addtocounter{theorem}{1}\begin{equation}} 
\def\eeq{\end{equation}} 

\def\beginproof{\noindent{\em Proof. }}
\def\Halmos{$\sqcup\llap{$\sqcap$}$}
\def\endproof{\unskip\nobreak\quad\Halmos\bigskip} 

\def\big{\bigskip}
 
\def\itt#1{{\em #1\/}}  

\def\a{\alpha}
\def\ad{\mbox{\rm ad}}
\def\Ad{\mbox{\rm Ad}}

\def\Ag{\mbox{\rm Aut}_k({\mathfrak g})}
\def\andd{\quad\mbox{and}\quad}
\def\Aut{\mbox{\rm Aut}}
\def\Aute{\mbox{\rm Aut}^{\mbox{\scriptsize e}}}
\def\Autz{\mbox{\rm Aut}^{0}}
\def\b{\mathfrak b}
\def\bl{\Big(}
\def\br{\Big)}
\def\bu{\mbox{\boldmath $u$}}
\def\bv{\mbox{\boldmath $v$}}
\def \C {\mathbb C}
\def\c{{\mathfrak c}}
\def\ck{^{\vee}} 
\def\Ctd{\mbox{\rm Ctd}}
\def\Del{\Delta}
\def\DelRe{\Delta^{\mbox{\scriptsize re}}}
\def\End{\mbox{\rm End}}
\def\ep{\varepsilon}
\def\f{{\mathfrak f}}
\def\fm{(\cdot,\cdot)}
\def\F{{\cal F}}
\def\g{{\mathfrak g}}
\def\G{\Gamma}
\def\gK{{\mathfrak g}(K)}
\def\GL{\mbox{\rm GL}}
\def\gR{{\mathfrak g}(R)}
\def\gS{{\mathfrak g}(S)} 

\def\gpK{{\mathfrak g}'(K)}
\def\h{{\mathfrak h}}
\def\Hom{\mbox{\rm Hom}}
\def\Hone{\mbox{\rm H}^1}
\def\i{{\bar \imath}}
\def\id{\mbox{\rm id}}
\def\implies{\Rightarrow}
\def\isom{\simeq}
\def\j{{\bar \jmath}}
\def\jfrak{{\mathfrak j}}
\def\k{k} 
\def\kunit{\k^\times}
\def\kalg{k\mbox{-alg}}
\def\L{{\cal L}}
\def\Lp{L}
\def\Lpg{L({\mathfrak g}, \sigma)}
\def\Lpgi{L({\mathfrak g}, \sigma^{-1})}
\def\Lpgp{L({\mathfrak g'}, \sigma |_{\g'})}

\def\M{{\cal M}}

\def\modcat{\mbox{-mod}}
\def\n{{\mathfrak n}}
\def\Out{\mbox{\rm Out}}
\def\pre#1{\, {}^{#1}}
\def\ot{\otimes}
\def\Sdm{{S_{dm}}}
\def\set#1{\{#1\}}
\def\sg{\sigma}
\def\suchthat{\mid}
\def\supth{{}^{\mbox{th}}}
\def\tah{\hat\tau}
\def\th{\theta}
\def\u{{\mathfrak u}}
\def\ui{\underline{i}}
\def\uo{\underline{1}}
\def\vp{\varphi}
\def\vpt{\tilde\varphi}

\def\Z{{\mathbb Z}}
\def\Zone{\mbox{\rm Z}^1}
\def\zt{\zeta}

\begin{document}
\title{Covering Algebras II: \\Isomorphism of Loop 
algebras}
\author{Bruce Allison, Stephen Berman and Arturo 
Pianzola\thanks{The authors gratefully acknowledge the support
of the Natural Sciences and Engineering Research Council of Canada.}}
\date{}
\maketitle
\begin{abstract} \noindent
This  paper studies the loop algebras that arise 
from  pairs consisting of a symmetrizable Kac-Moody Lie algebra $\g$ and a finite order automorphism $\sigma$ of $\g$. We obtain necessary and sufficient conditions for two such loop algebras to be isomorphic. 
\end{abstract}

\section{Introduction}

If $\sigma$ is a finite order automorphism of a Lie algebra $\g$,
Victor Kac \cite{K1} introduced the (twisted) loop algebra of $\g$ determined by $\sigma$.
Kac was interested in loop algebras since, in the case when
the base algebra $\g$ is finite dimensional simple, 
loop algebras provide explicit constructions (or realizations) of  affine Kac-Moody Lie algebras.  Our interest stems from the fact that when $\g$ is an affine Kac-Moody Lie algebra, loop algebras of $\g$ are used to provide realizations of extended affine Lie algebras of nullity 2 \cite{Wak,Po,A,SY}. Since the core of any extended affine Lie algebra is graded by a finite root system \cite{BGK,AG},
loop algebras based on affine algebras
also provide interesting examples of the root graded Lie algebras studied recently by several authors 
\cite{BM,BZ,ABG,BS}.

In this paper we solve the isomorphism problem for loop algebras of symmetrizable
Kac Moody Lie algebras.
That is, we derive necessary and sufficient conditions for  two such loop algebras to be isomorphic.  In order to accomplish this, we view  loop algebras as $S/R$-forms
of untwisted loop algebras, where $R$ is the ring of Laurent polynomials and $S/R$
is a Galois extension of rings.  We use this perspective in a number of ways---in particular,
it allows us to employ some of the machinery of nonabelian
Galois cohomology in our  arguments.

To be more precise about our results, let $\k$ be an algebraically closed field
of characteristic 0, let $\g$ be the Kac-Moody Lie algebra over $\k$
based on an indecomposable symmetrizable generalized Cartan matrix $A$,
and let $\sg$ be an automorphism of $\g$ of period $m$.
Let 
$$\g_{\i}=\{x \in \g \suchthat \sigma (x) = \zt_m^i x \},$$
where $\zt_m$ is a primitive $m^{\mbox{th}}$ root of unity in $\k$,
$i \in \Z$, and $i \to \i$ is the natural map of $\Z \to 
\Z_m$ (here $\Z_m$ denotes the integers modulo $m$). The loop algebra of 
$\g$ relative to $\sg$ is the subalgebra
$$\Lp(\g, \sigma) :=\bigoplus_{i \in \Z}\g_{\i} \ot_\k z^i$$
of the Lie algebra
$$\g(S) := \g\ot_\k S,$$
where $S := \k[z,z^{-1}]$ is the ring of Laurent
polynomials over $\k$.   As the notation suggests, the isomorphism class 
of $\Lpg$ does not depend on the choice of period $m$ (provided  that the primitive
roots $\zt_m$, $m\ge 1$, are compatible--see~(\ref{compat})).

In our first main result, Theorem \ref{Theorem:Invariance}, we show that $\g$ is an isomorphism invariant of
$\Lpg$.  Hence in the  isomorphism problem we may fix $\g$
and determine necessary and sufficient conditions for the isomorphism
of two loop algebras $\Lp(\g,\sg_1)$ and $\Lp(\g,\sg_2)$.

To discuss these conditions, let 
$\Aut(\g)$ denote the automorphism group 
of $\g$, let $\Aut(A)$ be the group of graph automorphisms of $A$ considered as a subgroup of $\Aut(\g)$,
and let $\omega$ be the Chevalley involution in $\Aut(\g)$. Define
$$\Out(A) := \cases{\Aut(A),& if $A$ has finite type\cr \langle \omega \rangle \times \Aut(A),&otherwise.\cr}$$
As usual $\Out(A)$ is called the outer automorphism group of $\g$ and there is a projection map
$$p: \Aut(\g) \rightarrow \Out(A).$$
(See \S \ref{Sec:Aut} for more details.)
Our second main result, Theorem \ref{Thm:Main}, says that if $\sg_1$ and $\sg_2$ are two automorphisms of $\g$ satisfying
$ \sg_1^m = \sg_2^m = \id$,  then 
$$\Lp(\g,\sg_1) \isom_k \Lp(\g,\sg_2) \iff p(\sg_1) \sim p(\sg_2),$$
where  $\isom_k$ denotes isomorphism over $\k$ and
$p(\sg_1) \sim p(\sg_2)$ means that
$p(\sg_1)$ is conjugate in $\Out(A)$ to either $p(\sg_1)$ or $p(\sg_2)^{-1}$. Notice that $\sg_1$ and $\sg_2$ need not have the same order but must only satisfy $\sg_1^m = \sg_2^m = \id$, and
so our result applies to any pair of finite order automorphisms (with $m$ chosen suitably).

In the  case when the base algebra $\g$ is finite dimensional, the above results
follow from work of Kac  on loop algebras \cite[Chapter 8]{K2}
and the Peterson-Kac conjugacy theorem for Cartan subalgebras of 
$\Lpg$ \cite[Theorem 2]{PK}.  In our more general context
we do not have conjugacy results in $\Lp(\g,\sg)$ to work with and so 
we must develop other tools. 
(The reader will notice however that conjugacy results in $\g$
are used frequently either directly or indirectly.) 
Thus, even when $\g$ is finite dimensional, our work gives
new proofs of some known results.  In particular we obtain
a proof, that does not use conjugacy in affine algebras,
of the nonisomorphism of the algebras that
appear in the classification of affine Kac-Moody Lie algebras   
(see Remark \ref{Rem:Affine}). 

To discuss our methods, we first mention a 
fact that is of central importance in this work. 
The loop algebra $\Lpg$
is an algebra over the ring $R=\k[t,t^{-1}]$, where $t = z^m$,
and furthermore the algebra $\Lpg\ot_R S$
is isomorphic as an algebra over $S$ to $\g(S)$.
Otherwise said, $\Lpg$ is an $S/R$-form of  $\g(R):= \g\ot_k R$.
(In fact the same is true when $\g$ is taken to be an arbitrary Lie algebra  over $\k$.)
This interpretation of loop algebras as forms
has an important consequence.  It allows us to 
view the $R$-isomorphism class of the loop algebra
$\Lpg$ as an element of the cohomology set $H^1(\Gamma, \Aut_S(\g (S)))$, 
where $\Gamma$ is the group of integers mod $m$ and $\Aut_S(\g (S))$ is the group of automorphisms of 
$\g (S)$ as an $S$-algebra. 
This  in turn allows us to use techniques from nonabelian  Galois cohomology to study isomorphism
over $R$ of two loop algebras of $\g$ (relative to automorphisms of period $m$).

Our interest though is in $\k$-isomorphism of loop algebras.  So it becomes
necessary to relate $\k$-isomorphism with $R$-isomorphism.
For this the key observation is that the centroid of $\Lp(\g',\sg|_{\g'})$ is canonically isomorphic to the algebra~$R$. This fact, together with a cohomological calculation, allows us to 
show that two loop algebras $\Lp(\g_1,\sg_1)$ and $\Lp(\g_2,\sg_2)$ are 
isomorphic as $\k$-algebras if and only if $\Lp(\g_1,\sg_1)$
is isomorphic as an $R$-algebra to either $\Lp(\g_2,\sg_2)$ or $\Lp(\g_2,\sg_2^{-1})$.

The tools just described are sufficient to prove the first main theorem
which states that  $\g$ is an isomorphism invariant of $\Lpg$.  They also
allow us to prove that the condition $p(\sg_1) \sim p(\sg_2)$
is necessary for isomorphism of  $\Lp(\g,\sg_1)$ and $\Lp(\g,\sg_2)$.

To prove the converse, additional tools, inspired by Kac's original argument
when $\g$ is finite dimensional \cite[Proposition 8.5]{K2}, are required.  
Among those tools are some results from the paper of Kac and Wang \cite{KW}
that allow us to  assume that the automorphisms of $\g$ that we are working with have a particularly nice Gantmacher like factorization. This in turn allows us to reduce 
our problem to considering the following situation. Let $\h$ be 
a fixed split Cartan subalgebra of $\g$ and let $\tau, \rho \in \Aut(\g)$. 
Assume that the  following three conditions hold: 
\begin{itemize}
\item[\rm (i)] $\tau^m=\rho^m = \id$, $\tau \rho=\rho \tau.$
\item[\rm (ii)] $\rho$ fixes $\h$ pointwise
\item[\rm (iii)] $\tau(\h)=\h.$
\end{itemize}
We are then able to prove what we think of as an 
``erasing result''. Namely we are able to show that in this situation we have 
$\Lp(\g, \tau \rho)  \isom_\k \Lp(\g, \tau)$.
Thus, we can erase $\rho$ and still stay in the same isomorphism class.
(The erasing result just indicated is actually a special case of a 
more general theorem that we prove in Section \ref{Sec:Erasing}.)
This then allows us to complete the proof that the condition $p(\sg_1) \sim p(\sg_2)$
is sufficient for isomorphism of  $\Lp(\g_1,\sg_1)$ and $\Lp(\g_2,\sg_2)$.

This present paper is the second in a sequence of three papers that
study loop algebras $\Lp(\g,\sg)$ and their ``double extensions'' 
$\mbox{\rm Aff}(\g,\sg) = \Lpg \oplus \k c \oplus \k d$, when $\g$ is an affine Kac-Moody Lie algebra and $\sg$ is a finite order automorphism
of~$\g$.  (We  refer to both $\Lpg$ and $\mbox{\rm Aff}(\g,\sg)$
as ``covering algebras'' of $\g$.)

In each of the first two papers we allow more general base algebras $\g$
than just affine algebras.
Indeed, our first paper \cite{A} assumed that $\g$ is a tame
extended affine Lie algebra (EALA for short)  
and determined when $\mbox{\rm Aff}(\g,\sg)$ is again a tame EALA.
Since affine Kac-Moody Lie algebras exactly correspond to  
tame EALAs of nullity 1, that paper applied in particular
to the case when $\g$ is affine.
In the present paper, we are not able to work in the general class
of tame EALAs since we need conjugacy results for Cartan 
subalgebras of $\g$ as well as the information concerning automorphisms which 
flows from this.
It would have been natural then, in this present paper, to simply assume that $\g$
is affine. However, since virtually no extra work results,
we have chosen to work 
in the more general setting of symmetrizable Kac-Moody Lie algebras. Thus, 
the first two papers are logically independent  and one can be read  
without knowing the other. It is in the third paper that we plan to put these results together.

The paper is organized so that the reader need only know some standard facts about symmetrizable Kac-Moody Lie  algebras. The deeper facts about automorphisms which we need are recalled, with references given, when  used. The same is true about the cohomological set up and  arguments. We have thus attempted to make this paper quite  self-contained.

\medskip
\noindent{\em  Assumptions and notation.}
\smallskip

We close this introduction by setting our basic assumptions and notation.

Although  our main interest is in algebras over algebraically closed
fields of characteristic 0, it is important (because we use base field extension
arguments) for us to relax this assumption somewhat at first.
We assume  that $\k$ \itt{is a field of characteristic 0 and that $\k$ contains
all roots of unity.}  The latter assumption means that for each integer
$m\ge 1$, $\k$ contains a primitive $m^{\mbox{th}}$ root of unity $\zt_m$.
We assume that these primitive roots of unity are  compatible in the sense that 
\beq
\label{compat}
\zt_{\ell m}^\ell = \zt_m
\eeq
for all $\ell, m\ge 1$.
It  is easy to see that such a choice is always possible, and we fix the sequence
$\{\zt_m\}_{m\ge 1}$ from now on.

We denote the group of units of $\k$ by $\kunit$.

Let $\kalg$ denote the category of unital commutative associative
$\k$-algebras.  If $K\in \kalg$ and $\L$ is a algebra over $K$, 
we use the  notation $\Aut_K(\L)$ for the group
of all $K$-algebra automorphisms of $\L$.   If $\M$
is a $K$-subalgebra of $\L$, then we set 
$$\Aut_K(\L;\M) : = \set{\tau\in \Aut_K(\L) \suchthat
\tau(x) = x \mbox{ for all } x\in \M}.$$
Also, if $\L_1$ and
$\L_2$ are algebras over $K$, we write $\L_1\isom_K \L_2$
to mean that $\L_1$ and $\L_2$ are isomorphic as algebras
over $K$.  If $\L$ is a Lie 
algebra over $K$  and $\M$ is a $K$-subalgebra of $\L$,
we denote the \itt{centralizer} of $\M$ in $\L$ by $C_{\L}(\M)$.
We use the notation $Z(\L) := C_\L(\L)$ for the \itt{centre} of $\L$.

We denote the tensor product $\ot_\k$ simply by $\ot$.
If $K\in \kalg$ and $V$ is a vector space over $\k$,
we use the functorial notation
$${V}(K) := {V}\ot K.$$
${V}(K)$ is a $K$-module with $K$-action given by 
$a(x\ot b) = x\ot ab$ for $x\in {V}$
and $a,b\in K$.   If $\g$ is an algebra over $\k$, then $\gK$
is an algebra over $K$.   Finally, if $\g$ is a Lie algebra over $\k$ and
$\f$ is a subalgebra of $\g$, then
$C_{\gK} (\f(K)) = (C_\g(\f))(K)$.  In particular, $Z(\gK) = (Z(\g))(K)$.

We fix a positive integer $m$. ($m$ will be a period
of the automorphisms investigated in this work.)
Let 
$$\G := \G_m := \Z /m\Z  = \set{\i : i\in \Z}$$
be the group of integers mod $m$, where $\i$ denotes the congruence class represented by $i$.

For  our purposes, the most important object in $\kalg$ will
be the algebra 
$$R = \k[t,t^{-1}]$$
of Laurent polynomials over $\k$ in the variable $t$.
We will also work extensively with the following algebra extension $S/R$.
Define
$$S := S_m := \k[z,z^{-1}]$$
to be the algebra of Laurent polynomials over $\k$ in the variable $z$,
and identify $R$ as a $\k$-subalgebra of $S$ via
$$t = z^m \quad\mbox{and hence}\quad t^{-1} = z^{-m}.$$
In this way $S$ is an $R$-algebra (which depends on $m$).

For $i\in \Z$, there is a unique $R$-algebra automorphism of $S$
which we denote by $\ui$ so that 
$$\ui(z) = \zt_m^i z.$$
Then the map $\i \mapsto \ui$ is an isomorphism of $\G$ onto  $\Aut_\k(S;R) = \Aut_R(S)$.
 
\goodbreak

\section{Loop algebras}

\label{Sec:Loop}
Suppose that $\g$ is a Lie algebra over $\k$ and
that $\sg$ is an automorphism of $\g$ of period $m$.  (The order of $\sg$ might not
be equal to $m$, but of course it is a divisor of~$m$.)
We begin by recalling the definition of the loop algebra of $\g$
relative to~$\sg$.

The algebra $\g$ can be 
decomposed into eigenspaces for $\sg$ as 
$$\g=  \bigoplus_{\i\in\G}\g_\i,$$
where 
$$\g_\i := \set{ x \in \g \suchthat \sg(x)= \zt_m^i x}.$$
Note that $\g_{\bar 0}$ is the subalgebra 
$\g^{\sg}$ of fixed points of $\sg$ in $\g$.

Now $\gS$ is a Lie algebra over $S$.  In  $\g(S)$ we set
$$
\Lp_m(\g,\sigma) =
\bigoplus_{i \in \Z}\g_{\i} \otimes   z^i.$$
Note  that $\Lp_m(\g,\sigma)$ is the set of fixed points in $\g(S)$ of the $R$-algebra
automorphism $\sg \ot \uo^{-1}$, and hence $\Lp_m(\g,\sigma)$ is an $R$-subalgebra of 
$\g(S)$.  Thus $\Lp_m(\g,\sigma)$ is an $R$-algebra and hence also a $\k$-algebra.
We call $\Lp_m(\g,\sigma)$ the \itt{loop algebra of $\g$
relative to~$\sg$}.

\begin{rem}  In this paper we consider loop algebras only in the class of Lie algebras.
However loop algebras can be defined in the same way in any class of algebras that allows base
ring extension.  (See \cite{P2}.)  Moreover the formal properties proved in Sections \ref{Sec:Loop}--\ref{Sec:Erasing} also hold in this generality (with the same proofs).
\end{rem}

\begin{rem}   The loop algebra just defined
depends on the choice of the primitive root of unity
$\zt_m$.  In fact, if we temporarily denote $\Lp_m(\g,\sg)$ by
$\Lp_m(\g,\sg,\zt_m)$,  we have
$\Lp_m(\g,\sg,\zt_m^a) = \Lp_m(\g,\sg^b,\zt_m)$ for $a,b\in \Z$ with $ab \equiv 1 \pmod m$.
Nevertheless, we will suppress the dependence of $\Lp_m(\g,\sg)$ on $\zt_m$ in our notation since we
have fixed the sequence $\{\zt_m\}_{m\ge 1}$ once and for all at the  outset.
\end{rem}

We now see that the $R$-isomorphism class of $\Lp_m(\g,\sigma)$ does not depend
on the choice of the period $m$ for~$\sg$.

\begin{lm}
\label{Lemma:Loopm} Suppose that $m'$ is another period of $\sg$.
Then
$$\Lp_{m'}(\g,\sigma) \isom_R \Lp_m(\g,\sigma) .$$
\end{lm}

\beginproof  We can assume without loss of generality that $m' = d$, where $d$ is a divisor of $m$.
Then $\zt_d = \zt_m^{e}$, where $e = \frac m d$.

Let $\varphi : \g \ot S_d \to \g \ot S_m$ be the $\k$-linear map so that 
$\varphi(x\otimes z^j) = x \otimes z^{ej}$
for $x\in \g$ and $j\in \Z$.  Then $\varphi$ is a $\k$-algebra monomorphism.  Moreover
\begin{eqnarray*}
&\varphi(t^i(x\ot z^j)) &= \varphi(x\ot t^i z^j) = \varphi(x\ot z^{di+j})=  x\ot z^{edi+ej}
\\&&=x\ot z^{mi+ej} = x\ot t^i z^{ej} = t^i \varphi(x\ot z^j)
\end{eqnarray*}
for $x\in \g$, $i,j\in \Z$. Hence $\varphi$ is $R$-linear.
Finally, since the eigenvalues of $\sg$ are contained in $\{\zt_d^0,\zt_d^1,\dots,\zt_d^{d-1}\} =
\{\zt_m^0,\zt_m^e,\dots,\zt_m^{e(d-1)}\}$, it is clear that $\varphi$ maps $\Lp_d(\g,\sigma)$
onto $\Lp_m(\g,\sigma)$.
\endproof

In view of this last lemma, there is no danger in abbreviating
our notation for loop algebras as follows:
$$\Lpg := \Lp_m(\g,\sigma).$$

We will need the following simple facts about loop algebras:

\begin{lm} 
\label{Lemma:LoopBasic}  With the above notation we have:
\begin{itemize} 
\item[\rm (a)] $\Lpgi \isom_k \Lpg$.
\item[\rm (b)] If $\tau\in\Aut_k(\g)$, then
$\Lp(\g, \tau\sg\tau^{-1}) \isom_R \Lpg$.
\end{itemize}
\end{lm}

\beginproof  (a): Let $\kappa : S \to S$ be the unique $\k$-algebra automorphism
of $S$ so that $\kappa(z) = z^{-1}$.  Then
$\id\ot\kappa$ is a 
$\k$-algebra automorphism of $\gS$ which maps $\Lpg$
onto $\Lpgi$.

(b): Let $\tau\in\Aut_k(\g)$.  Then $\tau\ot\id$ is an 
$R$-algebra automorphism of $\gS$ which maps $\Lpg$
onto $\Lp(\g, \tau\sg\tau^{-1})$.
\endproof

\section{Loop algebras as forms}
\label{Sec:Form}

Suppose that $\g$ is a Lie algebra over $\k$.
In this section we show how to interpret loop algebras
as $S/R$ forms of the $R$-algebra $\gR$.  
Since $S/R$ forms have a cohomological interpretation,
this gives us a cohomological tool to study loop algebras.

An  $S/R$ \itt{form} of $\gR$ is a Lie algebra $\F$ over $R$
such that ${\F}\ot_R S \isom \gS$ as Lie $S$-algebras.  
Since $S$ is a free $R$-module with basis 
$1,z,\dots,z^{m-1}$, the following lemma is clear.

\begin{lm}
\label{Lemma:Form}
If $\F$ is an $S/R$ form of $\gR$, then $\F$
can be identified with an $R$-subalgebra $\F$ of $\gS$ such
that each element of $\gS$ is uniquely expressible in the form
$\sum_{i=0}^{m-1} z^iy_i,$
where $y_i\in {\F}$ for $i=0,\dots,m-1$.  Conversely,
any $R$ subalgebra $\F$ of $\gS$ with this property is an $S/R$ form of $\gR$.
\end{lm}

Henceforth we will assume 
that all $S/R$ forms of $\gR$  occur as $R$-subalgebras
of $\gS$ with the property indicated in the lemma.  

Next we interpret $S/R$ forms
of $\gS$ cohomologically. For this purpose, we recall the definition
of  $\Hone(\G,F)$, when $F$ is a $\G$-group (see \cite[Section 5.1]{Ser}).

Suppose that $F$ is a group.  We say that $\G$-\itt{acts} on $F$ if
there is a map $(\i,f) \mapsto \pre{\i}f$ of $\G\times F \to F$
satisfying $\pre{\i}(\pre{\j}x) = \pre{\i+\j}x$, $\pre{\bar 0} x = x$ and 
$\pre{\i}(xy) = (\pre{\i}x) (\pre{\i}y)$ for $\i, \j\in \G$, $x,y\in F$.
In that case, $F$ is called a $\G$-\itt{group}. 

If $\G$ acts on a group $F$, 
a \itt{1-cocycle} on $\G$
in $F$ is a map $u : \G \to F$ satisfying the cocycle condition 
\beq
\label{Eq:Cocycle}
u_{\i+\j} = u_\i \pre{\i}u_\j
\eeq
for $\i,\j\in\G$, where the image of $\i$ under $u$ is denoted by $u_\i$.  The set of all such 1-cocycles is denoted by $\Zone(\G,F)$.  Two 1-cocycles $u$ and $v$ are said to be \itt{cohomologous} if there exists $f\in F$ so that
$$v_\i = f^{-1} u_\i \pre{\i}f$$
for $\i\in \G$. This is an equivalence relation on $\Zone(\G,F)$.  
If $u\in \Zone(\G,F)$ 
the equivalence class containing $u$
will be denoted by $\bu$ (the corresponding boldface character).  The set of all equivalence
classes (called cohomology classes) in $\Zone(\G,F)$ will be denoted by $\Hone(\G,F)$.  
$\Hone(\G,F)$ is a pointed set whose distinguished element is the class represented
by the constant 1-cocycle $\i \mapsto 1_F$.
Since $\G$ will be fixed
throughout the paper, we use the abbreviation
$$\Hone(F) := \Hone(\G,F).$$

If $\G$ acts on two groups $F_1$ and $F_2$,  a $\G$-\itt{morphism}  from $F_1$ to $F_2$
is a group homomorphism $\varphi : F_1 \to F_2$ that preserves the $\G$-action.
If $\varphi : F_1 \to F_2$ is $\G$-morphism, then
there is an induced  morphism
$\Hone(\varphi) : \Hone(F_1) \to \Hone(F_2)$ of pointed sets.  Indeed, $\Hone(\varphi)$ maps the class
represented by the 1-cocycle $u$ to the class represented by the 1- cocycle~$\varphi u$.

We define an action of $\G$ on $\Aut_S(\gS)$ by
\beq
\label{Eq:Action}
\pre{\i}\tau := (\id\ot\ui)\tau(\id\ot \ui)^{-1}
\eeq
for $\i\in\G$ and $\tau\in \Aut_S(\gS)$.
Using this action, we have:

\begin{prop} 
\label{Prop:H1} 
$R$-isomorphism classes of
$S/R$ forms of $\gR$ are in one-to-one correspondence with cohomology classes in  
$\Hone(\Aut_S(\gS)).$
Explicitly this correspondence is as follows:  If $u$ is a 1-cocycle on $\G$ in $\Aut_S(\gS)$, then the cohomology class $\bu$ corresponds to the $R$-isomorphism class of the $S/R$ 
form
\beq
\label{Eq:DefFrom}
\gS_u := \set{x\in \gS : (u_\i(\id\ot\ui))x = x \mbox{ for } \i\in\G}
\eeq
of $\gS$.
\end{prop}

\beginproof  
This proposition follows from general considerations using the
fact that $S/R$ is a Galois extension of rings with Galois group $\G$ (see \cite{Wat} or
\cite[\S III.2, Example 2.6]{M}). 
However, because of the special nature of the extension $S/R$, we can give 
a more elementary proof which we outline for the interest of the reader.
This is based on the approach to  $K/\k$ forms described in \cite[\S X.2]{J1} 
when $K/\k$ is a finite Galois extension of fields.

Suppose first that $u\in \Zone(\G,\Aut_S(\gS))$.
Let 
$\F = (\gS)_u$, where $(\gS)_u$ is defined by (\ref{Eq:DefFrom}).  We want to show that 
$\F$ is an $S/R$ form of $\gR$. To see this, set
$w_\i = u_\i(1\ot\ui)$ for $\i\in \G$.  Then $w_\i$ is an $R$-algebra automorphism of $\gS$
for $\i\in \G$.  Moreover, the cocycle condition (\ref{Eq:Cocycle}) translates to
$w_\i w_\j= w_{\i+\j}$ for $\i,\j\in\G$.  Thus $\set{w_\i \suchthat \i\in\G}$
is a cyclic group with generator $w_{\bar 1}$.  Consequently
$\F$ is the $R$-subalgebra of fixed points of $w_{\bar 1}$ in $\gS$.
Next, for $\i\in\G$, let 
$(\gS)_\i$ be the $\zt_m^i$-eigenspace of $w_{\bar 1}$ in $\gS$.
Since ${w_{\bar 1}}^m = \id$, we have $\gS = \bigoplus_{\i\in\G}(\gS)_\i$.
Furthermore, since $w_{\bar 1}$ is $\uo$-semilinear, we have $t^i (\gS)_{\bar 0} \subseteq
(\gS)_\i$ and $t^{-i} (\gS)_{\i} \subseteq (\gS)_{\bar 0}$ for $0\le i\le m-1$.
Thus $(\gS)_\i = t^i (\gS)_{\bar 0} = t^i\F$. 
So $\gS = \bigoplus_{\i\in\G}t^i\F$,
and hence, since $t^i$ is invertible in $S$,
$\F$ is an $S/R$ form of $\gR$.

Next suppose that $u,v$ are in $\Zone(\G,\Aut_S(\gS))$.  We claim that
$(\gS)_u$  and  $(\gS)_v$  are isomorphic as $R$-algebras if and only if $u$ and $v$ are cohomologous.
To see this, suppose first that $(\gS)_u$  and  $(\gS)_v$  are isomorphic as $R$-algebras.  Let $f:(\gS)_v\to(\gS)_u$ be an $R$-algebra isomorphism.  Since
$(\gS)_v$ and $(\gS)_u$ are $S/R$ forms of $\gR$, $f$ extends uniquely to an element $f\in \Aut_S(\gS)$. Then,
the $R$-algebra automorphisms 
$f^{-1}u_\i(1\ot\ui)f$
and 
$v_\i(1\ot \ui)$ of $\gS$ are both $\ui$-semilinear and both
fix the elements of $(\gS)_v$.  Hence, $ f^{-1}u_\i(1\ot\ui)f = v_\i(1\ot \ui)$ 
and so 
$v_\i = f^{-1} u_\i \pre{\i}f$ 
for $\i\in\G$.  Thus $u$ and $v$ are cohomologous. Conversely, suppose that $u$ and $v$ are cohomologous.  Let $f\in \Aut_S(\gS)$ with $v_\i = f^{-1} u_\i \pre{\i}f$ 
for $\i\in\G$.  Then $ f^{-1}u_\i(1\ot\ui)f = v_\i(1\ot \ui)$ for $\i\in\G$ and
so $f$ maps $(\gS)_v$ onto $(\gS)_u$.  Thus $(\gS)_u$ and $(\gS)_v$ are isomorphic
as $R$ algebras as claimed.  

So we have a well defined injective map from $\Hone(\Aut_S(\gS))$ into the set of $S/R$ forms of $\gR$
(regarded as $R$-subalgebras of $\gS$), namely the map that sends $\bu$ to the $R$-isomorphism
class of $(\gS)_u$.  It remains to show that this map is onto.  To do this suppose that $\F$ is
an $S/R$ form of $\gR$.  For $\i\in \G$, let $w_\i$ be the unique $R$-linear and $\i$-semilinear
map from $\gS$ to $\gS$ which fixes the elements of $\F$. Then $w_\i$ is an $R$-algebra automorphism of $\gS$ and $w_\i w_\j= w_{\i+\j}$ for $\i,\j\in\G$.  Set
$u_{\i} = w_\i(\id\ot \ui)^{-1}$
for $\i\in\G$.  Then $u\in \Zone(\G,\Aut_S(\gS))$ and one easily checks that
$\F = (\gS)_u$.
\endproof

The following theorem gives our description of loop algebras
as $S/R$ forms of $\gR$.

\begin{thm} 
\label{Thm:Loop} Suppose that $\g$ is a Lie algebra and that $\sg$
is an automorphism of $\g$ of period $m$.  Then
\begin{itemize}
\item[\rm (a)] $\Lpg$ is a free $R$-module of rank equal to $\dim_\k(\g)$.
\item[\rm (b)] $\Lpg$ is an $S/R$-form of $\gR$, and consequently
each element of $\gS$ can be written uniquely in the form
$\sum_{i=0}^{m-1} t^iy_i,$
where $y_i\in \Lpg$ for $i=0,\dots,m-1$.
\item[\rm (c)] The cohomology class in $\Hone(\Aut_S(\gS))$ that corresponds to the
$R$-iso\-mor\-phism class of $\Lpg$ is  the class represented by the $1$-cocycle
$\i \mapsto \sg^{-i}\ot \id$.
\end{itemize}
\end{thm}

\beginproof (a): If $B_\i$ is a $\k$-basis for $\g_\i$ for $\i\in \G$, then
$\cup_{i=0}^{m-1} B_\i\ot z^i$  is an $R$-basis for $\Lpg$.

(b) and (c): Set $u_\i = \sg^{-i}\ot \id$ for $\i\in \G$.  One checks easily
that $u$ is a 1-cocycle on $\G$ in $\Aut_S(\gS)$ and that $(\gS)_u = \Lpg$.
In view of Proposition \ref{Prop:H1} this proves both (b) and (c).
\endproof

We note that it is also easy to prove (b) in the theorem directly.

\section{$\k$-isomorphism versus $R$-isomorphism}
\label{Sec:Compare}

In the previous section, we have seen how to interpret $R$-isomorphism classes
of loop algebras cohomologically.  Since our real interest is in $\k$-isomorphism
of loop algebras, we need to understand the connection between $\k$-isomorphism
and $R$-isomorphism.  In this section, we explore this for loop algebras 
of Lie algebras that are perfect and central over $\k$.

Suppose $K\in \kalg$ and  $\L$ is a Lie algebra over $K$.  
The
\itt{derived algebra} $\L'=[\L,\L]$
of $\L$ is the $\Z$-span 
of the set of commutators of $\L$. $\L'$ is a $K$-subalgebra of $\L$.
We say that $\L$ is \itt{perfect}
if $\L' = \L$.  Note that we may regard the Lie algebra $\L$ over $K$
as a Lie algebra over $\k$.
This however does not change the derived algebra.  Hence, $\L$ is perfect
as a Lie algebra over $K$ if and only if it is perfect as Lie algebra over $\k$.

Again suppose $K\in \kalg$ and  $\L$ is a Lie algebra over $K$. We define 
$\lambda_\L: K \to 
\End_{K\modcat}(\L)$ by
$$\lambda_\L(a)(x) = ax$$
for $a\in K$ and $x\in\L$.  
We also define
$$\Ctd_K(\L) = \set{\chi \in \End_{K\modcat}(\L) \suchthat \chi([x,y]) = [x,\chi(y)] \mbox{ for all }
x,y \in \L}.$$
$\Ctd_K(\L)$ is a $K$-subalgebra of $\End_{K\modcat}(\L)$ 
called the \itt{centroid} of $\L$ over $K$. 
In general, we have $\lambda_\L(K)\subseteq \Ctd_K(\L)$.  We say that
$\L$ is \itt{central} over $K$ if $\Ctd_K(\L)=\lambda_\L(K)$. 

We next prove two lemmas about base ring restriction and extension.

\begin{lm} 
\label{Lemma:Centroid1}
Suppose that $K\in \kalg$ and $\L$ is a perfect Lie algebra over $K$.
Then $\Ctd_\k(\L) = \Ctd_K(\L)$.
\end{lm}
\beginproof  We must show that any element $\chi\in \Ctd_\k(\L)$ is $K$-linear.
Indeed for $a\in K$ and $x,y\in \L$, we have
$\chi(a[x,y]) = \chi([ax,y]) = [ax,\chi(y)] = a[x,\chi(y)] = a\chi([x,y])$.
\endproof

\begin{lm}
\label{Lemma:Centroid2}
Suppose that $\g$ is a Lie algebra over $\k$ and $K\in\kalg$.  
\begin{itemize}
\item[\rm (a)] $\gK' = \gpK$. 
\item[\rm (b)] If $\g$ is perfect, then $\gK$ is perfect.
\item[\rm (c)]  If $\g$ is central over $\k$, then $\gK$ is central over $K$.
\item[\rm (d)] If $\g$ is perfect and central over $\k$, then 
$$\Ctd_\k(\gK) = \Ctd_K(\gK) = \lambda_{\gK}(K).$$
\end{itemize}
\end{lm}

\beginproof
(a) is clear and (b) follows from (a).  To prove (c), we choose a basis $\{b_i\}_{i\in I}$ for the $\k$-space
$K$.  Let $\chi\in \Ctd_K(\gK)$.  Then  there exist unique elements 
$\chi_i\in \End_{\k\modcat}(\g)$,
$i\in I$, so that:
$$\chi(x\ot 1) = \sum_{i\in I} \chi_i(x)\ot b_i$$
for $x\in \g$, where 
for each $x\in \g$ only finitely many $\chi_i(x)$'s are nonzero.  One checks
easily that $\chi_i\in \Ctd_\k(\g)$ and so $\chi_i = \lambda_\g(a_i)$ for some $a_i\in \k$, $i\in I$.
Now we may assume that $\g$ is not zero and hence choose a nonzero $x_0\in \g$.
Then 
$\chi_i(x_0) = a_ix_0$ for $i\in I$.  Thus all but finitely many $a_i$'s are zero.
Set $a := \sum_{i\in I} a_ib_i \in K$.  Then, for $x\in \g$, we have
$\chi(x\ot 1) =  \sum_{i\in I} \chi_i(x)\ot b_i =  \sum_{i\in I} a_ix\ot b_i
= \sum_{i\in I} x\ot a_i b_i = x\ot a = \lambda_{\gK}(a)(x\ot 1)$.
Hence, $\chi= \lambda_{\gK}(a)$.  This proves (c).  (d) follows from (b),
(c) and Lemma \ref{Lemma:Centroid1}.\endproof

We are now ready to study $\k$-isomorphism versus $R$-isomorphism for loop algebras.
We first obtain an analog of the previous lemma for 
loop algebras.

\begin{lm}
\label{Lemma:CLoop}
Suppose that $\g$ is a Lie algebra over $\k$ and that $\sg$ is an automorphism
of period $m$ of $\g$.
\begin{itemize}
\item[\rm (a)] $\Lpg' = \Lpgp$. 
\item[\rm (b)] If $\g$ is perfect, then $\Lpg$ is perfect.
\item[\rm (c)]  If $\g$ is central over $\k$, then $\Lpg$ is central over $R$.
\item[\rm (d)] If $\g$ is perfect and central over $\k$, then 
$$\Ctd_\k(\Lpg) = \Ctd_R(\Lpg) = \lambda_{\Lpg}(R).$$
\end{itemize}
\end{lm}

\beginproof (a):  It is clear that $\Lpg' \subseteq \Lpgp$.  On the other hand,
for all $i\in \Z$, we have 
$\g_\i' \ot z^i = \sum_{j=0}^{m-1} [\g_\j,\g_{\i - \j}]\ot z^i = 
\sum_{j=0}^{m-1} [\g_\j\ot z^j,\g_{\i - \j}\ot z^{i-j}]\subseteq \Lpg'$.

(b): This follows from (a).

(c): Let $\L = \Lpg$.  Then, by Theorem \ref{Thm:Loop},
$\L$ is an $S/R$ form of $\gS$. 
We must show that $\Ctd_R(\L) \subseteq \lambda_\L(R)$.  For this purpose,
let $\chi\in \Ctd_R(\L)$.  Since $\L$ is an $S/R$ form of $\gS$, there exists
a unique $S$-linear map $\hat \chi : \gS \to \gS$ that extends $\chi$.
Explicitly, $\hat \chi$ is given by
$$\hat\chi(\sum_{i=0}^{m-1} z^iy_i) = \sum_{i=0}^{m-1} z^i\chi(y_i)$$
for $y_0,\dots,y_{m-1}\in \L$.  One checks easily that $\hat \chi\in \Ctd_S(\gS)$.
Thus, by Lemma \ref{Lemma:Centroid2}(c), $\hat\chi = \lambda_{\gS}(s)$ for some $s\in S$.
Write $s = \sum_{i=0}^{m-1}z^ir_i$, where $r_i\in R$.  Then, for $y\in \L$, we have
$\chi(y) = \hat\chi(y) = sy = \sum_{i=0}^{m-1}z^ir_iy$.  Hence,
$\chi(y) = r_0y$ and so $\chi = \lambda_\L(r_0)$.

(d): This follows from (b), (c) and Lemma \ref{Lemma:Centroid1}.
 \endproof

\begin{lm}
\label{Lemma:Induced}
Let  $\g_i$ be a nontrivial Lie algebra over $\k$ that is perfect and central over $\k$,
and let $\sg_i$ be an automorphism of $\g_i$ of period $m$, $i=1,2,3$.
\begin{itemize}
\item[\rm (a)]
If $\vp : \Lp(\g_1,\sg_1) \to \Lp(\g_2,\sg_2)$ is a $\k$-algebra isomorphism,
then there exists a unique $\vpt\in \Aut_\k(R)$ so that
$$\vp(rx) = \vpt(r)\vp(x)$$
for $x\in \Lp(\g_1,\sg_1)$ and $r\in R$.
\item[\rm (b)] If $\vp : \Lp(\g_1,\sg_1) \to \Lp(\g_2,\sg_2)$ and $\psi : \Lp(\g_2,\sg_2) \to \Lp(\g_3,\sg_3)$
are $\k$-algebra isomorphisms, then $\widetilde{(\psi\vp)} = \tilde \psi\tilde \vp $.
\end{itemize}
\end{lm}

\beginproof  (b) follows from (a), and so we only need  to prove (a).
 Let $\L_i = \Lp(\g_i,\sg_i)$, $i=1,2$.
Then the map 
\beq
\label{Eq:Induced}
\chi \mapsto \vp\chi\vp^{-1}
\eeq
is a $\k$-algebra isomorphism
from $\Ctd_k(\L_1)$ onto $\Ctd_k(\L_2)$.
So, by Lemma \ref{Lemma:CLoop}(d), the map (\ref{Eq:Induced}) is a $\k$-algebra
isomorphism of $\lambda_{\L_1}(R)$ onto $\lambda_{\L_2}(R)$.
But $\L_1$ and $\L_2$ are nontrivial free $R$-modules
by Theorem \ref{Thm:Loop}(a).  Thus the map $r\mapsto \lambda_{\L_i}(r)$
is a $\k$-algebra isomorphism from $R$ onto $\lambda_{\L_i}(R)$, $i=1,2$.
It follows that there exists a unique $\k$-algebra automorphism $\vpt$ of $R$ so that the following diagram commutes: 

\setlength{\unitlength}{.75pt}
\begin{picture}(120,100)(-150,0)
\put(0,0){\makebox(0,0){$\lambda_{\L_1}(R)$}}
\put(100,0){\makebox(0,0){$\lambda_{\L_2}(R)$}}
\put(0,65){\makebox(0,0){$R$}}
\put(100,65){\makebox(0,0){$R$}}
\put(25,0){\vector(1,0){50}}
\put(25,65){\vector(1,0){50}}
\put(0,50){\vector(0,-1){40}}
\put(100,50){\vector(0,-1){40}}
\put(-15,32){\makebox(0,0){$\lambda_{\L_1}$}}
\put(115,32){\makebox(0,0){$\lambda_{\L_2}$}}
\put(50,75){\makebox(0,0){$\vpt$}}
\end{picture}

\bigskip\noindent
Here the bottom map is the map (\ref{Eq:Induced}).
This completes the proof of (a).
\endproof

\begin{thm}
\label{Thm:Compare} Suppose that
$\k^m = \k$ (that  is every element of $\k$ is an $m\supth$  power  of an element of $\k$).
Let  $\g_i$ be a Lie algebra over $\k$ that is perfect
and central over $\k$, and let $\sg_i$ be an automorphism of period $m$
of $\g_i$, $i=1,2$.  Then $\Lp(\g_1,\sg_1)\isom_\k \Lp(\g_2,\sg_2)$ if and only if $\Lp(\g_1,\sg_1)\isom_R \Lp(\g_2,\sg_2)$ or $\Lp(\g_1,\sg_1)\isom_R \Lp(\g_2,\sg_2^{-1})$
\end{thm}

\beginproof  The sufficiency of these conditions follows from Lemma \ref{Lemma:LoopBasic}(a).  To prove
the necessity suppose that $\L_1\isom_\k \L_2$, where $\L_i := \Lp(\g_i,\sg_i)$ for $i=1,2$.
Let $\vp : \L_1 \to \L_2$ be an isomorphism of $\k$-algebras.  
We can assume that $\g_1$ and $\g_2$ are nontrivial, and  hence we can construct $\vpt\in \Aut_\k(R)$
as in Lemma \ref{Lemma:Induced}(a).  Then, since $\vpt$ is a $k$-automorphism of $R = \k[t,t^{-1}]$, 
$\vpt$ maps $t$ to a nonzero $\k$-multiple of $t$ or $t^{-1}$.

Let $\kappa : S \to S$ be the $\k$-algebra automorphism
of $S$ so that $\kappa(z) = z^{-1}$.  As we saw in the proof of Lemma \ref{Lemma:LoopBasic}(a),
$\psi := (\id\ot\kappa)|_{\L_2}$ is a 
$\k$-algebra isomorphism of  $\L_2$ onto $\Lp(\g_2,\sg_2^{-1})$.  Furthermore,
it is clear that $\tilde\psi = \kappa|_R$.  Hence, $\tilde\psi(t) = t^{-1}$.
But $\widetilde{(\psi\vp)} = \tilde \psi\tilde \vp$ by Lemma \ref{Lemma:Induced}(b).
Thus, replacing $\vp$ by $\psi\vp$ and $\sg_2$ by $\sg_2^{-1}$ if necessary, we can
assume that 
$$\vpt(t) = at,$$
where $a$ is a nonzero element of $\k$.

Choose $b\in\k$ so that $b^{-m} = a$.  Let $\ep$ be the  $\k$-algebra
automorphism of $S$ satisfying $\ep(z) = bz$.  Then $\tau := (\id\ot\ep)|_{\L_2}$ is a $\k$-algebra  automorphism of $\L_2$ and we have $\tilde\tau := \ep|_R$.
So $\tau\vp$ is a $\k$-algebra isomorphism of $\L_1$ onto $\L_2$
and $\widetilde{\tau\vp}(t) = \tilde\tau(\vpt(t)) = \tilde\tau(at) = \ep(at) =
\ep(az^m) = 
a(\ep(z))^m = a(bz)^m = z^m = t$.  Thus $\widetilde{\tau\vp}= \id$ and so $\tau\vp$
is an $R$-algebra isomorphism of $\L_1$ onto $\L_2$.
\endproof

\section{Erasing}
\label{Sec:Erasing}

Later  in this paper we will need to establish isomorphism of loop algebras
$\Lp(\g,\tau\rho)$ and $\Lp(\g,\tau)$ under suitable assumptions on $\g$,
$\tau$ and $\rho$.  In other words we will need to 
``erase'' the automorphism $\rho$.  In this section
we prove a rather general erasing theorem that will handle the cases
that we encounter.

If $\g = \oplus_{\a\in Q} \g_\a$ is a $Q$-graded Lie algebra where 
$Q$ is an abelian group and $r\in \Hom(Q,\k^\times)$,
we define $\Ad(r)\in \Aut_\k(\g)$ 
by
$$\Ad(r)(x_\a) := r(\a) x_\a$$
for $x_\a\in \g_\a$, $\a\in Q$.

We can now state our erasing theorem. The authors are grateful to  John Faulkner for conversations
about this theorem. He showed us the proof presented below which is 
significantly shorter 
than our original proof. 

\begin{thm}
\label{Thm:Erasing} {\rm (Erasing theorem)} Suppose that $\g = \oplus_{\a\in Q} \g_\a$ is a $Q$-graded Lie algebra where $Q$ is a free abelian group that is generated by the support
$\{\a \in Q : \g_\a \ne 0\}$ of $\g$. Suppose that $\tau, \rho$ are automorphisms of
$\g$  that satisfy the following assumptions:
\begin{itemize}
\item[\rm (i)] $\tau^m=\rho^m = \id$, $\tau \rho=\rho \tau.$
\item[\rm (ii)] $\rho = \Ad(r)$ for some $r\in \Hom(Q,\k^\times)$.
\item[\rm (iii)] There exists $\tah\in \Aut(Q)$ so that $\tau(\g_\a) = \g_{\tah(\a)}$
for $\a\in Q$.
\end{itemize}
Then
\beq
\label{looperase2}
\Lp(\g,\tau\rho) \isom_R \Lp(\g,\tau).
\eeq
\end{thm}

\beginproof Before beginning, notice that since $\tau^m = \id$ and since $Q$ is generated by the support of $\g$, it follows that $\tah^m$ has period $m$.  Let $d$ be any period
of $\tah$ (for example $d = m$).
We will prove that
\beq
\label{looperase}
\Lp_{dm}(\g,\tau\rho) \isom_R \Lp_{dm}(\g,\tau).
\eeq
By Lemma~\ref{Lemma:Loopm}, this implies (\ref{looperase2}). 

We prove (\ref{looperase})
by constructing an $\Sdm$-automorphism  
$\varphi$ of $\g \ot \Sdm$ with the property that
\beq
\label{conj}
\varphi(\tau\ot  \uo^{-1})\varphi^{-1} = \tau \rho \ot   \uo^{-1},
\eeq
where $\uo$ is the $\k$-automorphism of $\Sdm$ such that 
$$\uo(z) = \zt_{dm} z.$$  
Indeed, once we have constructed such an automorphism $\varphi$, we obtain
(\ref{looperase}) since
$\Lp_{dm}(\g,\tau)$ is the set of fixed points in
$\g\ot \Sdm$ of $\tau\ot \uo^{-1}$
and $\Lp_{dm}(\g,\tau\rho)$ is the set of fixed points in
$\g\ot \Sdm$ of $\tau\rho\ot \uo^{-1}$  (see Section \ref{Sec:Loop}).

So it remains to construct $\varphi\in \Aut_\Sdm(\g\ot \Sdm)$ satisfying
(\ref{conj}).  First of all, if $\a\in Q$ and $x_\a\in \g_\a$ we have
$x_\a = \rho^m(x_\a) = r(\a)^mx_\a$.  Thus $r(\a)\in \langle \zt_m \rangle$
for all $\a\in Q$ for which $\g_\a \ne 0$.  Since $Q$ is generated
by the support of $\g$, it follows that $r(Q)\subseteq \langle \zt_m \rangle$.
Consequently since $Q$ is free, there exists $a\in \Hom(Q,\Z)$ so that
$r(\a) = \zt_m^{a(\a)}$ for $\a\in Q$.  Hence
$$\rho(x_\a) =  r(\a) x_\a = \zt_m^{a(\a)} x_\a = \zt_{dm}^{da(\a)} x_\a$$
for $\a\in Q$ and $x_\a\in \g_\a$.

Now let $\varphi : \g\ot \Sdm \to \g\ot \Sdm$ be the $\k$-linear map so that
$$\varphi(x_\a \ot z^i) = \zt_{dm}^{b(\a)} x_\a\ot z^{i+c(\a)},$$
for $\a\in Q$ and $i\in \Z$, where $b,c\in \Hom(Q,\Z)$ are to be chosen 
below (depending on $a$ and $\tah$).  Since
$b$ and $c$ are homomorphisms and since $\g$ is $Q$-graded,
it follows that $\varphi\in \Aut_\Sdm(\g\ot \Sdm)$.  So we must
check (\ref{conj}) (with appropriately chosen $b$ and $c$).

Now for $x_\a\in \g_\a$ and $i\in \Z$, we have
\begin{eqnarray*}
&(\varphi(\tau\ot  \uo^{-1})\varphi^{-1})(x_\a\ot z^i) &=
\zt_{dm}^{-b(\a)}(\varphi(\tau\ot  \uo^{-1}))(x_\a\ot z^{i-c(\a)})
\\
&&=\zt_{dm}^{-b(\a)+c(\a)-i}\varphi(\tau(x_\a)\ot z^{i-c(\a)})
\\
&&
=\zt_{dm}^{-b(\a)+c(\a)-i+b(\tah(\a))}\tau(x_\a)\ot z^{i-c(\a)+c(\tah(\a))}
\end{eqnarray*}
and
$$(\tau\rho\ot  \uo^{-1})(x_\a\ot z^i) =\zt_{dm}^{-i}\tau\rho(x_\a)\ot z^i
=\zt_{dm}^{da(\a)-i}\tau(x_\a)\ot z^i.$$
So to insure that (\ref{conj}) holds it is sufficient  to choose
$b,c\in \Hom(Q,\Z)$ so that
\beq
\label{need1}
b\circ(\tah   - \id) + c = da
\eeq
and
\beq
\label{need2}
c\circ(\tah - \id) = 0,
\eeq
where $\circ$ denotes composition.

Next we define a polynomial over $\Z$ (in the indeterminant $x$) by
$$g(x) = \sum_{i=0}^{d-2}(d-i-1)x^i.$$
(Define $g(x) = 0$ if $d=1$.)
One checks directly that $g(x)(x-1)^2 = x^d-dx +d-1$.  Therefore since
$\tah^d = \id$ we have
\beq
\label{use}
g(\tah)\circ(\tah - \id)^2 = -d(\tah-\id).
\eeq
Let
$$b = -a\circ g(\tah)$$
and then let
$$c = da - b\circ(\tah-\id).$$
Then certainly (\ref{need1}) holds.  Moreover,
$$c\circ(\tah-\id) = da\circ(\tah-\id) - b\circ(\tah-\id)^2
= da\circ(\tah-\id) + a\circ g(\tah)\circ(\tah-\id)^2 
=0$$
by (\ref{use}). Thus we have (\ref{need2}), and the proof is  complete.
\endproof

\begin{rem} 
Equation  \ref{conj}  shows that the
cocycles determining $\Lp_{dm}(\g,\tau)$ and $\Lp_{dm}(\g,\tau\rho)$ 
are cohomologous.
It would be interesting to prove this using general cohomological  principles.
\end{rem}

\section{Symmetrizable Kac-Moody Lie algebras}
\label{Sec:KM}

For the rest of the paper  we will assume that 
$\g$ is an indecomposable symmetrizable Kac-Moody Lie algebra.
We begin this section  section 
by recalling the facts that we will need about these algebras.  The reader is referred to \cite[Chapters 1--4]{K2} or [MP, Chapter 4] for facts stated in these two sections without proof.  (In \cite{K2}, $\k$ is assumed to be the field of complex numbers.  However, that assumption is not needed for the results that we describe.)

We begin by establishing the notation that will be used for the
rest of the paper.

Let $A = (a_{ij})$ be an $n\times n$ \itt{symmetrizable} generalized Cartan matrix (GCM).   So, by definition, $a_{ii} = 2$ for all $i$, $a_{ij}$ is a nonpositive integer for $i\ne j$, and there exists an invertible diagonal matrix $D = \mbox{diag}(\ep_1,\dots,\ep_n)$ over $\k$ so that 
$D^{-1}A$ is symmetric.  We can (and do) choose $D$ with rational entries. We further
assume that $A$ is \itt{indecomposable}. Hence, the matrix $D$ is unique up to a 
nonzero constant factor.  

Choose a realization  $(\h,\Pi,\Pi\ck)$ of $A$.  Thus, by definition,
$\h$ is a vector space over $\k$ of dimension $n + \mbox{corank}(A)$, $\Pi = \set{\a_1,\dots,\a_n}$ is a linearly independent subset of the dual space $\h^*$ of $\h$, 
$\Pi\ck = \set{\a_1\ck,\dots,\a_n\ck}$ is a linearly independent subset of $\h$ and $\a_j(\a_i\ck) = a_{ij}$ for all $i,j$.  ($(\h,\Pi,\Pi\ck)$ is called a minimal realization of $A$ in [MP].)

Let 
$$\g = \g(A)$$
be the \itt{Kac-Moody Lie algebra} determined by $A$.  That is 
$\g$ is the Lie algebra over $\k$ generated by $\h$ and the symbols $\set{e_i}_{i=1}^n$
and $\set{f_i}_{i=1}^n$ subject to the relations $[\h,\h] = \set{0}$, 
$[e_i,f_j] = \delta_{ij} \a_i\ck$, 
$[h,e_i] = \a_i(h)e_i$, $[h,f_i] = -\a_i(h)f_i$ and 
$\ad(e_i)^{1-a_{ij}}(e_j)= \ad(f_i)^{1-a_{ij}}(f_j) =0 \ (i\ne j)$.

Let $\Del$ be the set of roots of $\h$ in $\g$.  Then we have the root space decomposition
$\g = \bigoplus_{\a\in\Del\cup\set{0}}\g_\a$, where $\g_0 = \h$.  Let $\DelRe$ denote the set of
real roots in~$\Del$.  Let $Q = \sum_{i=1}^n \Z \a_i$ be the \itt{root lattice} of $\g$.

Let $\h'  := \h\cap \g' = \sum_{i=1}^n \k\a_i\ck$.  (This is an exception to our notation in Section~\ref{Sec:Compare}.
$\h'$  is of course not the derived algebra of $\h$.)
In that case we have $\g' = \h' \oplus(\bigoplus_{\a\in\Del}\g_\a)$.

Let  $\c$ be the centre of $\g'$.  Then, $\c$ is also the centralizer of  $\g'$ in $\g$,
and we have $\c = \{h\in \h : \a_i(h) = 0 \mbox{ for } i=1,\dots,n\}$.

Next let $\Aut(A)$ be the \itt{group of automorphisms of the GCM} $A$.  Thus,
by definition $\Aut(A)$ is the group of permutations $\nu$ of $\set{1,\dots,n}$ so that
$a_{\nu i,\nu j} = a_{i,j}$ for $1\le i,j\le n$.  It is important for our purposes
to identify $\Aut(A)$ as a subgroup of $\Ag$.  Indeed, it is shown in
\cite[\S 4.19]{KW} that there exists a subspace $\h''$ of $\h$ 
and an action of $\Aut(A)$ on $\h$
so that 
$$\h = \h'\oplus \h'',$$
$\nu(\h') = \h'$ and $\nu(\h'') = \h''$ for $\nu\in \Aut(A)$, and
$$ \nu(\a_i\ck) = \a_{\nu i}\ck \andd \a_{\nu i}(\nu(h'')) = \a_{i}(h'')$$
for $\nu\in \Aut(A)$, $1\le i\le n$ and $h''\in \h''$.
(The choice of $\h''$ is not unique.  However, once such a choice is made, the action of 
$\Aut(A)$ on $\h$ is determined.) 
One then has
$$\a_{\nu i}(\nu(h)) = \a_{i}(h)$$
for $\nu\in \Aut(A)$, $1\le i\le n$ and $h\in \h$.
Furthermore, the action of $\Aut(A)$ on $\h$ extends uniquely to an action of
$\Aut(A)$ on $\g$ by automorphisms so that 
$$\nu(e_i) = e_{\nu i} \andd \nu(f_i) = f_{\nu i},\quad$$
for $\nu\in \Aut(A)$ and $1\le i\le n$.  From now on we fix a choice of $\h''$
as above and hence regard $\Aut(A)$ as a subgroup of $\Aut_k(\g)$.

Since $A$ is symmetrizable, we can use 
the complementary 
space $\h''$ for $\h'$ in $\h$ to define an invariant
form $\fm$ on $\g$.  Indeed, there exists a unique invariant symmetric $\k$-bilinear form 
$\fm$ on $\g$ so that $(\h'',\h'') = \set{0}$ and
$$(\a_i\ck,h) = \ep_i \a_i(h)$$
for $1\le i \le n$ and $h\in \h$.    This form is nondegenerate on $\g$.
Moreover,  the automorphisms of $\g$ in $\Aut(A)$
preserve this  form (see  the argument in \cite[4.15]{KW}).

We will need two lemmas about the Lie algebra $\g$.  The first deals with
basefield extension.

\begin{lm}  
\label{Lemma:Extend}
Suppose that $K$ is an extension field of $\k$.  Then (with the obvious identifications) $(\h(K),\Pi \ot K,\Pi\ck\ot K)$ is a realization of $A$ over $K$,
and $\gK$ is isomorphic to the Kac-Moody Lie algebra over $K$ determined by $A$ using the 
realization $(\h(K),\Pi \ot K,\Pi\ck\ot K)$.
\end{lm} 
\beginproof  This is proved in \cite[\S 4.3]{MP}.  Since we need the explicit isomorphism,
we briefly describe the proof.
The first statement is clear.  To prove the second statement, we (temporarily)
let $\tilde\g$ be the Kac-Moody Lie algebra over $K$ determined by $A$ using the 
realization $(\h(K),\Pi \ot K,\Pi\ck\ot K)$.  Then $\tilde\g$ is generated by
$\h(K)$,  $\set{\tilde e_i}_{i=1}^n$
and $\set{\tilde f_i}_{i=1}^n$ with relations as above.  It follows from the universal
property of $\g$ and the tensor product
that there is a $K$-algebra homomorphism from $\gK$ to $\tilde g$ 
that restricts to the identity on $\h(K)$ and maps $e_i\ot 1$ to $\tilde e_i$ and
$f_i\ot 1$ to $\tilde f_i$.  One obtains an inverse of this homomorphism using
the universal property of $\tilde g$.
\endproof

The second lemma will be useful since it allows us to apply the results of Section
\ref{Sec:Compare} to the Lie algebra $\g'$.

\begin{lm}
\label{Lemma:GpPerfect}
The Lie algebra $\g'$ is perfect and central over $\k$.
\end{lm}

\beginproof  It is clear  (and well known) that $\g'$ is perfect since 
$\g'$ is generated by $\bigcup_{i=1}^n\{e_i,f_i\}$
and since $e_i = \frac 12 [\a_i\ck,e_i]$ and $f_i = -\frac 12 [\a_i\ck,f_i]$.

To show that $\g'$ is central over $\k$, suppose that $\chi\in \Ctd_\k(\g')$.  We first show that
\beq
\label{Eq:Chi}
\chi(\h') \subseteq \h'.
\eeq
For this purpose, let $g\in \h'$ and write 
$\chi(g) = x_0 +\sum_{\a\in \Del}x_\a$,
where $x_0\in \h'$ and $x_\a\in \g_\a$ for $\a\in \Del$.  
Suppose for contradiction that $x_\beta \ne 0$ for some $\beta\in \Del$.
Then $0 = \chi([g,g]) = [g,\chi(g)] = [g,x_0 +\sum_{\a\in \Del}x_\a]= 
\sum_{\a\in \Del}\a(g)x_\a$ and so $\beta(g) = 0$. Thus,
for $y_\beta\in \g_{-\beta}$, we have
$$[y_\beta,x_0 +\sum_{\a\in \Del}x_\a] = [y_\beta,\chi(g)] = \chi([y_\beta,g]) = 
\beta(g)\chi(y_\beta) = 0.$$
Hence, $[y_\beta,x_\beta] = 0$ for all $y_\beta\in \g_{-\beta}$.  This contradicts the assumption that $x_\beta \ne 0$ (see  \cite[Theorem 2.2(e)]{K1})
and so we have (\ref{Eq:Chi}).  Now fix  $i\in\set{1,\dots,n}$.
By (\ref{Eq:Chi}), $\chi(\a_i\ck)\in\h'$. Thus $$\chi(e_i) = \frac 12\chi([\a_i\ck,e_i]) = \frac 12[\chi(\a_i\ck),e_i] = \frac 12\a_i(\chi(\a_i\ck))e_i = ae_i,$$
where $a =  \frac 12\a_i(\chi(\a_i\ck))\in\k$.  But since $A$ is indecomposable,
the ideal of $\g'$ generated by $e_i$  is  $\g'$, and so $\chi(x) = ax$ for all $x\in\g'$
as desired.
\endproof

We  can now prove our first main theorem
which shows that $\g$ is an isomorphism  invariant of $\Lp(\g,\sg)$ and of 
$\Lp(\g',\sg)$.

\begin{thm}  
\label{Theorem:Invariance}
Suppose that $\k$ is a field of characteristic 0 that contains all roots of unity and
suppose that $\k^m = \k$. Let $\g_i$ be an indecomposable symmetrizable Kac-Moody Lie algebra for $i=1,2$.
\begin{itemize}
\item[\rm (a)] If $\sg_i$ is an automorphism of period $m$ of $\g_i$ for $i = 1,2$ and if
$\Lp(\g_1,\sg_1) \isom_k \Lp(\g_2,\sg_2)$ then $\g_1 \isom_k \g_2$. 
\item[\rm (b)] If $\sg_i$ is an automorphism of period $m$ of $\g_i'$ for $i = 1,2$
and if $\Lp(\g'_1,\sg_1) \isom_k \Lp(\g'_2,\sg_2)$ then $\g_1 \isom_k \g_2$. 
\end{itemize}
\end{thm}

\beginproof By Lemma \ref{Lemma:CLoop}(a), it is sufficient to prove (b).
Suppose that $\sg_i$ is an automorphism of period $m$ of $\g'$ for $i = 1,2$,
and suppose that $\Lp(\g'_1,\sg_1) \isom_k \Lp(\g'_2,\sg_2)$. 
By Lemma \ref{Lemma:GpPerfect} and Theorem \ref{Thm:Compare}, we have
$\Lp(\g'_1,\sg_1)\isom_R \Lp(\g'_2,\sg_2)$ or $\Lp(\g'_1,\sg_1)\isom_R \Lp(\g'_2,\sg_2^{-1})$.
We can assume that  $\Lp(\g'_1,\sg_1)\isom_R \Lp(\g'_2,\sg_2)$.
Hence $\Lp(\g'_1,\sg_1)\ot_R S \isom_S \Lp(\g'_2,\sg_2)\ot_R S$ 
and so, by Theorem \ref{Thm:Loop}(b), we have $\g'_1(S) \isom_S \g'_2(S)$.
Let $L$ be the quotient field of the integral domain $S$.
Then $\g'_1(S)\ot_S L \isom_S \g'_2(S)\ot_S L$ and so $\g'_1(L) \isom_L \g'_2(L)$.
Thus, by Lemma \ref{Lemma:Centroid2}(a), we have $(\g_1(L))' \isom_L (\g_2(L))'$.
Now by assumption, $\g_i$ is the Kac-Moody Lie algebra over $\k$
determined by an indecomposable symmetrizable
GCM $A_i$, $i=1,2$.   Hence, by Lemma \ref{Lemma:Extend},  $\g_i(L)$
is the Kac-Moody Lie algebra over $L$ determined by $A_i$, $i=1,2$.
Thus the fact that $(\g_1(L))' \isom_L (\g_2(L))'$ implies that
$A_1$ and  $A_2$ are the same up to a bijection of the index sets \cite[Theorem 2(b)]{PK}.
Hence $\g_1 \isom_\k \g_2$.\endproof

\begin{rem} If $\g_1$ and $\g_2$ are finite dimensional,
Theorem \ref{Theorem:Invariance} follows from Exercise 8.7 of \cite{K1} or from
the Peterson-Kac conjugacy theorem 
for split Cartan subalgebras of $\Lp(\g_2,\sg_2)$. 
\end{rem}

In view of Theorem \ref{Theorem:Invariance}(a), to study the question
of isomorphism of loop algebras of indecomposable symmetrizable Kac-Moody
algebras it suffices to consider the case when the base algebra $\g$
is fixed.  We do this for the rest of the  paper.

\section{The automorphism group of $\g$}
\label{Sec:Aut}

We  continue to assume that $\g = \g(A)$ is the Kac-Moody Lie algebra 
determined by an indecomposable symmetrizable GCM $A$.
We  next discuss the structure of the group $\Ag$.  We follow [KW] in this discussion.
We begin by introducing some subgroups of $\Ag$.

Recall first that $\Aut_\k(\g;\h) = \set{\tau\in\Ag \suchthat \tau(h) = h \mbox{ for all } h\in\h}$.  One can give an explicit description of the elements of
$\Aut_\k(\g;\h)$.  Let $\Hom(Q,\kunit)$ be the set of group homomorphisms
of the root lattice $Q$ into the group $\kunit$ of units  of $k$.
$\Hom(Q,\kunit)$ is a group under pointwise   multiplication.
If $r \in \Hom(Q,\kunit)$, we define $\Ad(r)\in \Aut_\k(\g;\h)$ (as in Section \ref{Sec:Erasing})
by
$$\Ad(r)(x) = r(\a)x$$
for $x\in \g_a$, $\a\in \Del\cup\set{0}$.  Then $\Ad : \Hom(Q,\kunit) \to
\Aut_\k(\g;\h)$ is an isomorphism of groups.

Recall next that $\Aut_\k(\g;\g') = \set{\tau\in\Ag \suchthat \tau(x) = x \mbox{ for all } x\in\g'}$.  Again one can give an explicit description of the elements of
$\Aut_\k(\g;\g')$.   Indeed, we have 
$$\g = \g'\oplus \h''.$$  Thus, if $\eta\in  \Hom_\k(\h'',\c)$, we may define a $\k$-linear map $\th_\eta : \g \to \g$ by 
$$\th_\eta(x + h) = x + h + \eta(h)$$
for $x\in\g'$ and $h\in \h''$.  Then the map $\eta\mapsto \th_\eta$ is a group isomorphism
of $\Hom_\k(\h'',\c)$ onto $\Aut_\k(\g;\g')$.  (We prove a slightly more general fact in 
Proposition \ref{Prop:Surjective}(a) below.)

Next let $\Aute_\k(\g)$ denote the subgroup of
$\Ag$ that is generated by $\set{\exp(\ad(x)) \suchthat x\in\g_\a,\ \a\in\DelRe}$.
Now $\Aut_\k(\g;\h)$ normalizes $\Aute_\k(\g)$, while $\Aut_\k(\g;\g')$ centralizes
both $\Aute_\k(\g)$ and $\Aut_\k(\g;\h)$.  Hence,
$$\Autz_\k(\g) := \Aute_\k(\g) \Aut_\k(\g;\h) \Aut_\k(\g;\g')$$
is a subgroup of $\Aut_k(\g)$.  We call $\Autz_\k(\g)$ the \itt{inner automorphism group of $\g$}.

Finally  let $\omega$ be the unique $\k$-automorphism
of $\g$  so that
$$\omega(e_i) = -f_i,\quad \omega(f_i) = -e_i\andd \omega(h) = -h$$
for $1\le i \le n$ and $h\in \h$.  Then $\omega$ is an element of order $2$ in $\Ag$
called the \itt{Chevalley automorphism} of $\g$.  Note that $\omega$ commutes with the
elements of $\Aut(A)$ (regarded as automorphisms of $\g$) and so we may set
$$\Out(A) := \cases{\Aut(A),& if $A$ has finite type\cr \langle \omega \rangle \times \Aut(A),&otherwise.\cr}$$
$\Out(A)$ is a finite subgroup of $\Ag$ that we call the \itt{outer automorphism group of~$\g$}.  
Observe that the automorphisms
of $\g$ in $\Out(A)$ preserve the form $\fm$ (since $\omega$ clearly does).
We will see below that $\Out(A)$ is a complement to  $\Autz_\k(\g)$ in $\Ag$.

We will also be interested in the groups $\Aut_\k(\g')$ and $\Aut_\k(\g'/\c)$.
We have group homomorphisms
\beq 
\label{map1}
\Ag \mapsto \Aut_\k(\g')
\eeq
and
\beq 
\label{map2}
\Aut_\k(\g') \mapsto \Aut_\k(\g'/\c),
\eeq
where (\ref{map1}) is the restriction map and (\ref{map2}) maps $\tau \in\Aut_\k(\g')$ to the automorphism of $\g'/\c$ induced by $\tau$. The homomorphism (\ref{map1}) has kernel $\Aut_\k(\g;\g')$ and, since $\g'$ is perfect, the homomorphism (\ref{map2}) is injective.

Let $\Autz_\k(\g')$ denote the image of $\Autz_\k(\g)$ under the map (\ref{map1}), and let
$\Autz_\k(\g'/\c)$ denote the image of $\Autz_\k(\g)$ under the composite of the maps
(\ref{map1}) and  (\ref{map2}). 

The map (\ref{map1}) restricted to $\Out(A)$
is an injection, and so, abusing notation, we denote the image of
$\Out(A)$ in both $\Aut_\k(\g')$ and  $\Aut_\k(\g'/\c)$  by $\Out(A)$.

If $H$ and $K$ are subgroups of a group $G$, we write (as usual) $G = H \rtimes K$
to mean that $H$ is a normal subgroup of $G$, $G = HK$ and $H\cap K = \{1\}$.

The following result on the structure of $\Aut_\k(\g)$ is due to Peterson and 
Kac~\cite{PK}.

\begin{prop}
\label{Prop:Aut} 
We have $\Aut_\k(\g'/\c)  = \Autz_\k(\g'/\c) \rtimes \Out(A)$, $ \Aut_\k(\g')= \Autz_\k(\g') \rtimes \Out(A)$ 
and
\beq
\label{decomp}
\Ag = \Autz_\k(\g) \rtimes \Out(A).
\eeq
Consequently the homomorphism (\ref{map1}) is surjective
and the homomorphism (\ref{map2}) is an isomorphism
\end{prop}

\beginproof
The decomposition $\Aut_\k(\g'/\c)  = \Autz_\k(\g'/\c) \rtimes \Out(A)$ is proved
in \cite[Theorem 2(c)]{PK} as a consequence of the conjugacy theorem for split Cartan
subalgebras of $\g'/\c$. (More details are given in \cite[1.25--1.32]{KW}).
The rest of the proposition follows easily from this.
\endproof

It follows from the fact that (\ref{map1}) is a surjective that we have
the exact sequence 
$$\set{1}\to \Aut_\k(\g;\g')\to \Ag \to \Aut_\k(\g')\to \set{1}.$$
Moreover, as we've noted, the group $\Aut_\k(\g;\g')$ in this sequence is isomorphic
to $\Hom_\k(\h'',\c)$.  In order to compute cohomology we will need
a  version of these facts  with the base ring extended.  We present this next.

If $K\in\kalg$, then we have
$\gK = \gpK\oplus \h''(K)$.  Thus, if $\eta$ is an element of  $\Hom_K(\h''(K),\c(K))$, we may define a $K$-linear map $\th_\eta : \gK \to \gK$ by 
$$\th_\eta(x + h) = x + h + \eta(h)$$
for $x\in\gpK$ and $h\in \h''(K)$.  

\begin{prop}  Let $K\in \kalg$.
\label{Prop:Surjective}
\begin{itemize}
\item[\rm (a)] The map $\eta \mapsto \th_\eta$ is an isomorphism from 
the group $\Hom_K(\h''(K),\c(K))$ onto the group $\Aut_K(\gK;\gpK)$.  
\item[\rm (b)] Suppose further that $K$ is an integral domain.
Then the restriction map $\Aut_K(\gK)\to \Aut_K(\gpK)$ is surjective,  and hence we have
the exact sequence
$$\set{1}\to \Aut_K(\gK;\gpK) \to \Aut_K(\gK) \to \Aut_K(\gpK)\to \set{1}.$$
\end{itemize}
\end{prop}

\beginproof (a):  If $\eta\in\Hom_K(\h''(K),\c(K))$, one easily checks that
$\th_\eta$ is a $K$-algebra homomorphism which fixes $\g'(K)$ pointwise.
Also, using the fact that $\c(K)$ is contained in $\gpK$, it is easy to check
that $$\th_{\eta_1}\th_{\eta_2}= \th_{\eta_1+\eta_2}$$
for $\eta_1,\eta_2\in\Hom_K(\h''(K),\c(K))$.
Therefore $\th_\eta\in\Aut_K(\gK;\gpK)$ for $\eta\in\Hom_K(\h''(K),\c(K))$,
and the map $\eta \mapsto \th_\eta$ is an injective group homomorphism from 
$\Hom_K(\h''(K),\c(K))$ into  $\Aut_K(\gK;\gpK)$.   It remains to
show that the image of this map is $\Aut_K(\gK;\gpK)$.
For this suppose that $\tau$ is an element of $\Aut_K(\gK;\gpK)$.
Then for $h\in \h''(K)$ and $x\in \gpK$, we have
$[\tau(h)-h,x] = [\tau(h),\tau(x)]-[h,x]=\tau([h,x])-[h,x] = 0$.
Hence, if $h\in \h''(K)$, $\tau(h)-h$ centralizes $\gpK$ .
But $\c$ is the centralizer of $\g'$ in $\g$, and so $\c(K)$
is the centralizer of $\gpK$ in $\gK$.  Thus $\tau(h)-h\in \c(K)$
for $h\in \h''(K)$.  Hence we may define $\eta:\h''(K)\to \c(K)$ by
$\eta(h) = \tau(h)-h$, in which case we have $\tau = \th_\eta$.

(b):  Suppose that $K$ is an integral domain and let $F$ be the quotient field
of $K$.  In this argument, we work in the Lie algebra $\g(F) = \g\ot F$ over the field
$F$.  Note that $\g(F)$ contains $\g(K)$ as a $K$-subalgebra.

We have $\c = \set{h\in \h : \a_i(h) =0 \mbox{ for } i=1,\dots,n}$.  Hence we can
choose a basis $g_1,\dots,g_n,c_1,\dots,c_p$ for $\h$ so that
$\a_i(g_j) = \delta_{ij}$ for $1\le i,j\le n$ and $c_1,\dots,c_p$ is a basis for 
$\c$.  Consider the $\k$-space
$$\u = \sum_{i=1}^n\k g_i + \sum_{\a\in\Del}\g_\a.$$
Then,
$$\g = \u\oplus\c.$$

Let $\n = \set{x\in \g(F) \suchthat [x,\gpK] \subseteq \gpK}$ be the normalizer
of $\gpK$ in $\g(F)$.  We first prove that
\beq 
\label{Eq:Surjective1}
\n = \u(K)\oplus \c(F).
\eeq
Indeed, the inclusion ``$\supseteq$'' is clear.  So we must show the inclusion ``$\subseteq$''.
Now $\gpK$ is stable under $\ad(\g\ot 1)$ and hence so is $\n$.  Thus $\n = \bigoplus_{\a\in\Del\cup\set{0}} \n\cap(\g_\a(F))$ and so it is enough to show that
$$\n\cap(\g_\a(F)) \subseteq \u(K)\oplus \c(F)$$
for $\a\in\Del\cup\set{0}$.  To show this let $x\in\n\cap(\g_\a(F))$, where $\a\in\Del\cup\set{0}$.  Suppose first that $\a=0$.  Then $x\in \n\cap(\h(F))$.  Hence,
$x = \sum_{i=1}^n g_i\ot a_i + \sum_{i=1}^p c_i\ot b_i$, where $a_i,b_i\in F$. 
Subtracting $\sum_{i=1}^p c_i\ot b_i$, we may assume that $x = \sum_{i=1}^n g_i\ot a_i$.
Then $[x,e_j\ot 1]\in \g'(K)$.  But 
$$[x,e_j\ot 1] = \sum_{i=1}^n [g_i,e_j]\ot a_i= \sum_{i=1}^n \a_j(g_i)e_j\ot a_i= e_j\ot a_j$$
and so $a_j\in K$.  Hence, $x\in \u(K)$.  Suppose next that $\a\in\Del$.  Choose
a basis $\set{x_i}_{i=1}^q$ for $\g_\a$ and a basis $\set{y_i}_{i=1}^q$ for $\g_{-\a}$
so that $(x_i,y_j) = \delta_{ij}$.  Then $x = \sum_{i=1}^q x_i\ot a_i$, where $a_i\in F$.
Now $[x,y_j\ot 1] \in \g'(K)$ and
$$[x,y_j\ot 1] = \sum_{i=1}^q [x_i,y_j]\ot a_i = [x_j,y_j]\ot a_j$$
for $1\le j\le q$.  Since $[x_j,y_j]$ is a nonzero element of $\h$ it follows that
$a_j\in K$ for $1\le j\le q$.  Thus $x\in \u(K)$ and we have proved (\ref{Eq:Surjective1}).

We are now ready to prove that the restriction map 
\beq
\label{Eq:Surjective2}
\Aut_K(\gK)\to \Aut_K(\gpK)
\eeq
is surjective.  To do this we identify $\Aut_K(\gK)$ as a subgroup of
$\Aut_F(\g(F))$, by identifying each element of $\Aut_K(\gK)$ with its unique extension
to an $F$-linear endomorphism of $\g(F)$.  Similarly, we identify $\Aut_K(\gpK)$ as a subgroup of
$\Aut_F(\g'(F))$.  Note that by Lemma \ref{Lemma:Extend},
$\g(F)$ is isomorphic to the Kac-Moody Lie algebra over the
field $F$ determined by the GCM $A$.  Hence, by Proposition \ref{Prop:Aut},
the restriction map
$\Aut_F(\g(F))\to \Aut_F((\g(F))')$ is surjective.  Thus, 
by Lemma \ref{Lemma:Centroid2}(a), the map $\Aut_F(\g(F))\to \Aut_F(\g'(F))$ is surjective.

To prove that the map (\ref{Eq:Surjective2}) is surjective, let $\rho\in \Aut_K(\gpK)$.
We may regard $\rho$ as an element of $\Aut_F(\g'(F))$, and hence 
there exists 
$\tau\in \Aut_F(\g(F))$ so that
\beq
\label{Eq:Surjective3}
\tau|_{\g'(F)} = \rho.
\eeq
Note that it follows that
\beq
\label{Eq:Surjective4}
\tau(\g'(K)) = \g'(K).
\eeq
In the remainder of the proof, we will show that $\tau\in \Aut_F(\g(F))$ can be chosen satisfying
(\ref{Eq:Surjective3}) (and hence (\ref{Eq:Surjective4})) as well as the condition
$$
\tau(\gK) = \gK.
$$
This will complete the proof since then $\tau$ (restricted to $\g(K)$) will be the 
required extension of $\rho$.  

Now, if $h\in \h''(K)$, we have
$$[\tau(h),\gpK] = [\tau(h), \tau(\gpK)] = \tau([h, \gpK])\subseteq \tau(\gpK) = \gpK.$$ 
Thus $\tau(\h''(K))\subseteq \n$ and so, by (\ref{Eq:Surjective1}), 
\beq
\label{Eq:Surjective5}
\tau(\h''(K)) \subseteq \u(K) \oplus \c(F).
\eeq
Consequently, we have
$$\tau(\h''(F)) \subseteq \u(F) \oplus \c(F),$$
and so there exist unique elements 
$$\gamma\in \Hom_F(\h''(F),\u(F)) \andd\ep\in \Hom_F(\h''(F),\c(F))$$
so that
$$\tau(h) = \gamma(h) + \ep(h)$$
for $h\in \h''$.
Moreover, by (\ref{Eq:Surjective5}), we have
\beq
\label{Eq:Surjective6}
\gamma(\h''(K)) \subseteq \u(K).
\eeq
Now $\c(F)$ is the centre of $\g(F)$ and so $\tau(\c(F)) = \c(F)$.  Thus we can
define $\eta\in \Hom_F(\h''(F),\c(F))$ by $\eta(h) = - \tau^{-1}\ep(h)$
for $h\in \h''(F)$.  We then can define $\th_\eta\in\Aut_F(\g(F);\g'(F))$ by
$\th_\eta(x+h) = x+h+\eta(h)$ for $x\in \g'(F)$, $h\in \h''(F)$.
In that case we have
$$\tau\th_\eta(h) = \tau(h+\eta(h)) = \tau(h) + \tau\eta(h)=\tau(h)-\ep(h)=\gamma(h)$$
for $h\in \h''(F)$.  Thus, by (\ref{Eq:Surjective6}), $\tau\th_\eta(\h''(K))\subseteq \u(K)$.
On the other hand, $\tau\th_\eta(\gpK) = \tau(\gpK) = \gpK$ by
(\ref{Eq:Surjective4}).  Hence, $\tau\th_\eta(\gK) \subseteq \gK$.
Thus, replacing $\tau$ by $\tau\th_\eta$, we have shown that we can choose 
$\tau\in \Aut_F(\g(F))$ so that $\tau|_{\g'(F)} = \rho$ and
$$\tau(\gK) \subseteq \gK.$$

Similarly, we can choose $\mu\in \Aut_F(\g(F))$ so that 
$$\mu|_{\g'(F)} = \rho^{-1} \andd
\mu(\gK) \subseteq \gK.$$
Set $\xi = \tau\mu\in \Aut_F(\g(F))$.  Then,
$\xi \in \Aut_F(\g(F);\g'(F))$ and 
\beq
\label{Eq:Surjective7}
\xi(\g(K)) \subseteq \g(K).
\eeq
Since $\xi \in \Aut_F(\g(F);\g'(F))$, it follows from (a) that there exists
an element $\eta\in \Hom_F(\h''(F),\c(F))$ so that 
$$\xi(h)= h+\eta(h)$$
for $h\in \h''(F)$.
If $h\in\h''(K)$, we have $\eta(h) = \xi(h)-h\in \g(K)$
and so $\eta(h)\in \g(K)\cap\c(F) = \c(K)$.  Hence,
$$\eta(\h''(K))\subseteq \c(K).$$
Thus,
\begin{eqnarray*}
&\gK &\subseteq\gpK+  \h''(K)
\subseteq \gpK + \xi(\h''(K)) + \eta(\h''(K))\\
&&\subseteq \gpK + \xi(\h''(K)) + \c(K)
\subseteq \gpK + \xi(\h''(K))\\
&&\subseteq\xi(\gpK)+ \xi(\h''(K))
\subseteq\xi(\gK)\\
\end{eqnarray*}
Combining this with (\ref{Eq:Surjective7}), we have $\xi(\g(K)) = \g(K)$.
Since $\mu$ and $\tau$ stabilize $\gK$ and $\xi = \tau\mu$, it follows that
$\tau(\g(K)) = \g(K)$ as desired.\endproof

\section{A necessary condition for isomorphism of loop algebras}
\label{Sec:Necessary}

Suppose again that $\g = \g(A)$ is the Kac-Moody Lie algebra 
determined by an indecomposable symmetrizable GCM $A$.
Our main goal in this section is to obtain a simple necessary condition for isomorphism
of loop algebras of $\g$ relative to finite period automorphisms of $\g$.

We will need the following lemma:

\begin{lm}
\label{Lemma:Key}
Suppose that $\sg_1$ and $\sg_2$
are automorphism of $\g$ of period $m$. If 
$\Lp(\g',\sg_1|_{\g'})\isom_R \Lp(\g',\sg_2|_{\g'})$,
then $\Lp(\g,\sg_1)\isom_R \Lp(\g,\sg_2)$.
\end{lm}

\beginproof  First of all, $S$ is an integral domain.  Hence, by Proposition \ref{Prop:Surjective}(b),
we have the exact sequence
\beq
\label{Eq:Key1}
\set{1}\to \Aut_S(\g(S);\g'(S)) \to \Aut_S(\g(S)) \stackrel{\pi} \to \Aut_S(\g'(S))\to \set{1},
\eeq
where  $\pi$ is the restriction map.
Also, $\G$ acts on $\Aut_S(\g(S))$ and $\Aut_S(\g'(S))$ by
$$
\pre{\i}\tau := (\id\ot\ui)\tau(\id\ot \ui)^{-1}.
$$
This action of $\G$ on $\Aut_S(\g(S))$ stabilizes $\Aut_S(\g(S);\g'(S))$. 
Hence,  the sequence (\ref{Eq:Key1}) is an exact sequence of $\G$-groups.

Assume now that $\Lp(\g',\sg_1|_{\g'})\isom_R \Lp(\g',\sg_2|_{\g'})$.
We define $1$-cocycles $u$ and $v$ on $\G$ in $\Aut_S(\g(S))$
by $u_\i = \sg_1^{-i}\ot \id$ and $v_\i = \sg_2^{-i}\ot \id$. 
Let $\bu$ and $\bv$ be the cohomology classes in $\Hone(\Aut_S(\g(S)))$ represented 
by $u$ and $v$ respectively.  Then, by Proposition \ref{Prop:H1} and Theorem \ref{Thm:Loop}(c), we have $\Hone(\pi)(\bu) = \Hone(\pi)(\bv)$.  We wish to
show that $\bu = \bv$.  Hence, it suffices to show: 
\beq
\label{Eq:Key2}
\mbox{The 
preimage under
$\Hone(\pi)$ of $\Hone(\pi)(\bu)$ is $\set{\bu}$. }
\eeq

There is a well known method that can be used to show (\ref{Eq:Key2}). This involves
``twisting'' the action of $\G$ by the cocycle $u$.
We define a new action ``$\cdot$'' of $\G$ on $\Aut_S(\g(S))$ by setting
$$
\i \cdot \tau := u_\i  \pre{\i}\tau u_\i^{-1} 
= ( \sg_1^{-i}\ot\ui)\tau(\sg_1^i\ot (\ui^{-1})),
$$
for $\i\in\G$ and $\tau\in \Aut_S(\g(S))$.  The action ``$\cdot$'' 
stabilizes $\Aut_S(\g(S);\g'(S))$.  We denote the $\G$-group
$\Aut_S(\g(S);\g'(S))$ with the action  ``$\cdot$''  by ${}_u\Aut_S(\g(S);\g'(S))$
Then, since the sequence (\ref{Eq:Key1}) is exact, 
it is known that (\ref{Eq:Key2})
follows from
\beq
\label{Eq:Key3}
\Hone({}_u\Aut_S(\g(S);\g'(S))) = \set{1}.
\eeq
(See \cite[Corollary 2, \S I5.5]{Ser}.)
Thus it is sufficient to prove (\ref{Eq:Key3}).

For the sake of brevity, we set $M := \Hom_S(\h''(S),\c(S))$.
Then, by Proposition \ref{Prop:Surjective}(a), 
the map $\eta \mapsto \th_\eta$ is an isomorphism from 
the group $M$ onto the group $\Aut_S(\g(S);\g'(S))$.  
We use this isomorphism to transport the action ``$\cdot$'' from
$\Aut_S(\g(S);\g'(S))$ to obtain an action, also denoted by ``$\cdot$'',
of $\G$ on $M$.
Thus we have 
$$\th_{\i\cdot \eta} = \i\cdot\th_\eta$$
for
$\i\in\G$ and $\eta\in M$.
We use the notation 
${}_uM$ for the additive group
$M$ with $\G$-action ``$\cdot$'' obtained in this way.
Then it is sufficient to show that 
$\Hone({}_uM) = \set{0}$.

Now there exists a unique element $\beta$  of period $m$ in 
$\GL_\k(\h'')$ and a unique element $\gamma$ of  
$\Hom_\k(\h'',\g')$ so that 
$$\sigma_1(h) = \beta(h) + \gamma(h)$$
for $h\in \h''$.  Moreover, 
$$\sigma_1^{-1}(h) = \beta^{-1}(h) + \delta(h)$$
for $h\in \h''$, where $\delta\in \Hom_\k(\h'',\g')$ is defined by
$\delta = -(\sigma_1^{-1}|_{\g'}) \gamma \beta^{-1}$.
Then, if $h\in\h''$ and $\eta\in M$,
we have
\begin{eqnarray*}
&\th_{\bar 1\cdot \eta}(h\ot 1)
&=(\bar 1\cdot \th_\eta)(h\ot 1)\\
&&=\bl( \sg_1^{-1}\ot\uo)\th_\eta(\sg_1\ot (\uo^{-1}))\br\bl h\ot 1 \br\\
&&= 
\bl( \sg_1^{-1}\ot\uo)\th_\eta\br
\bl\sg_1(h)\ot 1\br\\
&&= 
\bl( \sg_1^{-1}\ot\uo)\th_\eta\br 
\bl\beta(h)\ot 1+\gamma(h)\ot 1\br\\
&&= 
\bl \sg_1^{-1}\ot\uo\br
\bl\beta(h)\ot 1 +\eta(\beta(h)\ot 1) + \gamma(h)\ot 1\br\\
&&=h\ot 1 + (\delta\beta)(h)\ot 1\\
&&\qquad\qquad\qquad+ 
\bl(\sg_1^{-1}\ot \uo)\eta(\beta\ot \id)\br
\bl h\ot 1\br
+(\sg_1^{-1} \gamma)(h)\ot 1\\
&&=h\ot 1 +  
\bl(\sg_1^{-1}\ot \uo)\eta(\beta\ot \id)\br
\bl h\ot 1\br\\
&&= \th_{\tilde \eta}(h\ot 1),
\end{eqnarray*}
where $\tilde\eta \in  M$
is defined by $\tilde \eta = (\sg_1^{-1}|_{\c}\ot \uo)\eta(\beta\ot \id).$
Hence, we have $\bar 1\cdot \eta = (\sg_1^{-1}|_{\c}\ot \uo)\eta(\beta\ot \id)$, and so
\beq
\label{Eq:Key4}
\i \cdot \eta = (\sg_1^{-1}|_{\c}\ot \uo)^i\, \eta\, (\beta\ot \id)^i
\eeq
for $\eta\in M$ and $\i\in\G$. 

Now $_u{}M=  \Hom_S(\h''(S),\c(S))$ as a group and hence 
$_u{}M$ has the natural structure of a $\k$-vector space.
Moreover, if follows from (\ref{Eq:Key4}) 
that the $\G$-action  on ${}_uM$ is $\k$-linear.  
Hence, $\Hone({}_uM)$ is a $\k$-vector space.  
But since $\G$ is
a finite group of order $m$, we have 
$m\Hone({}_uM)=\{0\}$, by \cite[Theorem 6.14]{J2}, and hence
$\Hone({}_uM)=\{0\}$ as desired.
\endproof

\begin{thm}
\label{Thm:Equivalent}  Suppose that
$\k^m = \k$. Suppose  $\sg_1$ and $\sg_2$
are automorphisms of $\g$ of period $m$.
Then the following statements are equivalent:
\begin{itemize}
\item[\rm (a)] $\Lp(\g,\sg_1)\isom_\k \Lp(\g,\sg_2)$
\item[\rm (b)] $\Lp(\g',\sg_1|_{\g'})\isom_\k \Lp(\g',\sg_2|_{\g'})$
\item[\rm (c)] $\Lp(\g',\sg_1|_{\g'})\isom_R \Lp(\g',\sg_2|_{\g'})$ 
or $\Lp(\g',\sg_1|_{\g'})\isom_R \Lp(\g',\sg_2^{-1}|_{\g'})$
\item[\rm (d)] $\Lp(\g,\sg_1)\isom_R \Lp(\g,\sg_2)$ or $\Lp(\g,\sg_1)\isom_R \Lp(\g,\sg_2^{-1})$
\end{itemize}
\end{thm}
\beginproof
``(a) $\implies$ (b)''  follows from
Lemma \ref{Lemma:CLoop}(a). ``(b) $\implies$ (c)'' follows from Lemma \ref{Lemma:GpPerfect} and
Theorem \ref{Thm:Compare}. ``(c) $\implies$ (d)'' follows from Lemma \ref{Lemma:Key}.
Finally, ``(d) $\implies$ (a)'' follows from Lemma \ref{Lemma:LoopBasic}(a).
\endproof

Recall from (\ref{decomp}) that we have the decomposition
$$\Ag = \Autz_\k(\g) \rtimes \Out(A)$$
of $\Ag$.  Let $p_\k : \Ag \to \Out(A)$
be the projection
onto the second factor relative to this decomposition.  Often  we will drop
the subscript $\k$ from this notation and write 
\beq 
\label{project}
p = p_\k.
\eeq

We  next define an equivalence relation $\sim$ on $\Out(A)$ by defining:
\beq
\label{equiv}
\mu_1\sim\mu_2 \iff \mbox{$\mu_1$ is conjugate to either $\mu_2$ or $\mu_2^{-1}$ in $\Out(A)$}
\eeq
for  $\mu_1,\mu_2\in \Out(A)$.

We can now
give a simple necessary condition for isomorphism
of loop algebras of $\g$.  We will see in the final section that this condition is also sufficient.

\begin{thm}
\label{Thm:Necessary}  Suppose that
$\k^m = \k$.  Suppose that  $\sg_1$ and $\sg_2$
are automorphisms of $\g$ of period $m$. If $\Lp(\g,\sg_1)\isom_\k \Lp(\g,\sg_2)$,
then $p(\sg_1)\sim p(\sg_2)$.
\end{thm}

\beginproof  To prepare for the proof  we need some notation and some remarks.

First of all, to avoid confusion below, we denote the subgroup $\Out(A)$ of $\Ag$  by $\Out_\k(A)$. So we have
$\Ag = \Autz_\k(\g) \rtimes \Out_\k(A)$,
and $p_k$ is the projection onto the second factor.

Next let $L$ be the quotient field of  $S = \k[z,z^{-1}]$.
We identify $\Aut_S(\g(S))$ as a subgroup of $\Aut_L(\g(L))$ by identifying
each element of $\Aut_S(\g(S))$ with its unique extension
to an $L$-linear endomorphism of $\g(L)$.  Let $\ep : \Aut_S(\g(S))\to \Aut_L(\g(L))$
be the inclusion map obtained in this way.

If $i\in\Z$, the automorphism $\ui$ of $S$ satisfying $\ui(z) = \zt_m^i z$ extends
uniquely to an automorphism, which we also denote by $\ui$, of $L$.
This allows us to extend the action (\ref{Eq:Action}) of $\G$ on $\Aut_S(\gS)$ to an action of $\G$ on $\Aut_L(\g(L))$
by defining 
\beq
\label{Eq:Necessary1}
\pre{\i}\tau := (\id\ot\ui)\tau(\id\ot \ui)^{-1}.
\eeq
for $\i\in\G$ and $\tau\in \Aut_L(\g(L))$.  Then $\ep$ preserves
the action of $\G$.

Next, we use the isomorphism described in the proof of Lemma \ref{Lemma:Extend}
(with $K = L$) to identify $\g(L)$ with the Kac-Moody Lie algebra over $L$ determined
by the GCM $A$.  For the complement of $\h'(L)$ in $\h(L)$, we choose the space $\h''(L)$.
Then, by Proposition \ref{Prop:Aut}, we have a decomposition
$$
\Aut_L(\g(L)) = \Autz_L(\g(L)) \rtimes \Out_L(A).
$$
We let $p_L : \Aut_L(\g(L)) \to \Out_L(A)$ be the projection onto the
second factor in this decomposition. 

Note that
$$
\Out_L(A) = \Out_k(A) \ot \id.
$$
It is clear from this equation and (\ref{Eq:Necessary1}) that
$\G$ acts trivially on the subgroup $\Out_L(A)$ of $\Aut_L(\g(L))$.
Also from the definition of 
$\Autz_L(\g(L))$, it is clear that
$\Autz_L(\g(L))$ is stabilized by the action of $\G$.
Thus $p_L$ preserves the action of $\G$.

Finally note that if $\tau\in \Ag$, then $\tau\ot \id \in \Aut_L(\g(L)$ and
\beq
\label{Eq:Necessary2}
p_L(\tau\ot \id) = p_\k(\tau)\ot \id.
\eeq

We are now ready to give the proof of the theorem.
Suppose  that $\Lp(\g,\sg_1)\isom_\k \Lp(\g,\sg_2)$.  Then,
by Theorem \ref{Thm:Equivalent},
we have $\Lp(\g,\sg_1)\isom_R \Lp(\g,\sg_2)$ or $\Lp(\g,\sg_1)\isom_R \Lp(\g,\sg_2^{-1})$.
Replacing $\sg_2$ by $\sg_2^{-1}$ if necessary, we can assume that
\beq
\label{Eq:Necessary3}
\Lp(\g,\sg_1)\isom_R \Lp(\g,\sg_2).
\eeq
Let $u\in \Zone(\Aut_S(\g(S))$ and $v\in \Zone(\Aut_S(\g(S))$ be defined by
$u_\i = \sg_1^{-i} \ot \id$ and $v_\i = \sg_2^{-i} \ot \id$, $\i\in\G$.
Let $\bu,\bv\in  \Hone(\Aut_S(\g(S))$ be the cohomology classes represented by $u$
and $v$ respectively.  Then, by Proposition \ref{Prop:H1}, Theorem \ref{Thm:Loop}(c) and (\ref{Eq:Necessary3}), we have $\bu = \bv$.  But
$p_L\ep$ is a $\G$-homomorphism from $\Aut_S(\g(S))$ into $\Out_L(A)$.
Hence, we have the induced map $\Hone(p_L\ep) : \Hone(\Aut_S(\g(S))) \to \Hone(\Out_L(A))$.  Consequently, $\Hone(p_L\ep)(\bu) = \Hone(p_L\ep)(\bv)$.
But 
$$(p_L\ep)(u_\i) = p_L(\ep(\sg_1^{-i}\ot \id))= p_L(\sg_1^{-i}\ot \id) = p_\k(\sg_1^{-i})\ot \id$$
for $\i\in\G$, by (\ref{Eq:Necessary2}).  Thus,
$\Hone(p_L\ep)(\bu)$ is represented by the 1-cocycle $\i \mapsto  p_\k(\sg_1^{-i})\ot \id$.
Similarly, $\Hone(p_L\ep)(\bv)$ is represented by the 1-cocycle $\i \mapsto  p_\k(\sg_2^{-i})\ot \id$.  Since the action of $\G$ on $\Out_L(A) = \Out_\k(A)\ot \id$ is trivial, it follows 
that there exists $\mu\in \Out_\k(A)$ so that
$$(\mu \ot \id) (p_\k(\sg_1^{-i})\ot \id) (\mu\ot 1)^{-1} = p_\k(\sg_2^{-i})\ot \id$$
for $i\in \Z$.  Taking $i = -1$, we see that $p_\k(\sg_1)$ is conjugate
to $p_\k(\sg_2)$ in $\Out_\k(A)$.
\endproof

\begin{rem}
\label{Rem:Affine}
Suppose that $\k = \C$. Suppose that we let $\g = \g(A)$ run over all isomorphism classes
of finite dimensional simple Lie algebras, and for each $\g$ we let
$\sg$ run over a set of representatives of the $\sim$-equivalence classes
in $\Out(A) = \Aut(A)$.  For each of the pairs $(\g,\sg)$
we form  the loop algebra $\Lpg$.  
Kac  \cite[Chapter 8]{K2} has shown that this procedure yields 
the derived algebra modulo its centre 
of the affine algebras of types 
$$\mbox{A}_1^{(1)},\dots,\mbox{E}_8^{(1)},
\mbox{A}_2^{(2)},\dots,\mbox{E}_6^{(2)},\mbox{D}_4^{(3)}.
$$
That is, it yields the derived algebra modulo its centre of
all affine Kac-Moody Lie algebras up to isomorphism.
Now, it follows from Theorem \ref{Theorem:Invariance}(a) and  Theorem~\ref{Thm:Necessary} that
distinct pairs $(\g_1,\sg_1)$ and $(\g_2,\sg_2)$ yield
nonisomorphic loop algebras $\Lp(\g_1,\sg_1)$ and $\Lp(\g_2,\sg_2)$.
Consequently, we have a proof of the nonisomorphism of the algebras that 
appear in the classification of affine algebras. 
Of course this nonisomorphism result also follows from the Peterson--Kac conjugacy theorem for split Cartan subalgebras in affine algebras, but we are not using that deep result here  (when
the base algebra $\g$ is finite dimensional).  
\end{rem}

\section{Isomorphism of loop algebras}
\label{Sec:Converse}

In  this section we prove that the necessary condition for isomorphism of loop algebras given
in Theorem \ref{Thm:Necessary} is also sufficient.

We assume again that $\g = \g(A)$ is the Kac-Moody Lie algebra over $\k$
determined by an indecomposable symmetrizable GCM $A$,
and we use the notation of the previous section.
In particular, the projection $p : \Ag \to \Out(A)$ is defined by (\ref{project}),
and the equivalence relation $\sim$ on $\Out(A)$ is defined by (\ref{equiv}).

First of all, the following result follows immediately from the general erasing Theorem
\ref{Thm:Erasing}:

\begin{thm}
\label{Thm:ErasingKM}  Suppose that $\tau,\rho\in \Ag$ with
$\tau^m = \rho^m = \id$ and $\tau\rho = \rho\tau$.  Suppose also that 
$\rho\in \Aut_\k(\g;\h)$ and $\tau(\h) = \h$.  Then $\Lp(\g,\tau\rho) \isom \Lp(\g,\tau)$.
\end{thm}

If $\sg\in \Ag$, then, by \cite[Theorem 3]{PK}, $\sg(\b^+)$
is $\Aut_\k^e(\g)$-conjugate to either $\b^+$ or $\b^-$,
where $\b^+$ (resp.~$\b^-$) is the subalgebra of $\g$ generated
by $\{e_i\}_{i=1}^n$ (resp.~$\{f_i\}_{i=1}^n$).  The automorphism $\sg$ is
said to be of the {\it first kind\/} or the {\it second kind\/} accordingly.

The following result is proved using Theorem \ref{Thm:ErasingKM} along with detailed information from \cite{KW} about semisimple
automorphisms of the first and second kind of $\g$.  The results in \cite{KW} assume that $\k$ is algebraically closed
and so we assume that here.

\begin{thm} 
\label{Thm:Project}
Suppose that $\k$ is algebraically closed.
Let $\sg\in \Ag  $ with $\sg^m = \id$.
Then
$$\Lp(\g,\sg) \isom_R \Lp(\g,p(\sg)).$$
\end{thm}

\beginproof Suppost first that $\g$ is of the first kind.  Then  $\sg$, being of finite order,
is a semisimple automorphism of $\g$ of the first kind, and so, by \cite[Lemma 4.31]{KW},
$\sg$ is conjugate in $\Ag$ to an automorphism of $\g$ of the form $\nu\rho$, where
$\nu\in \Aut(A)$, $\rho\in \Aut_\k(\g;\h)$ and $\nu\rho= \rho\nu$.
By Lemma \ref{Lemma:LoopBasic}(b), we may assume that $\sg = \nu\rho$.  Then,
$p(\sg) = p(\nu)p(\rho) = \nu$.  Consequently, $\nu ^m = \id$ and so
$\rho^m = \id$.  Thus, by Theorem \ref{Thm:ErasingKM}, we have $\Lp(\g,\sg) \isom_R \Lp(\g,\nu)$,
and so $\Lp(\g,\sg) \isom_R \Lp(\g,p(\sg))$.

Before considering the case when $\sg$ is of the second kind, we need some notation.
Suppose for that purpose that $X$ is a subset of $\Pi$ of finite type (we allow $X$ to be
empty).   Let $\g_X$ be the
subalgebra of $\g$ generated by $\{e_i\}_{\a_i\in X}$ and $\{f_i\}_{\a_i\in X}$, and let $\h_X = \sum_{\a_i\in X} \k \a_i\ck$.  Then $\g_X$ is a finite dimensional semisimple subalgebra of $\g$
with Cartan subalgebra $\h_X$.  We have
$$\h = \h_X \oplus \h_X^\perp,$$ where $\h_X^\perp$ is the orthogonal complement of $\h_X$ in $\h$
(relative to the form $\fm$ defined in Section \ref{Sec:KM}).  Let
$E_X$ be the subgroup of $\Ag$ generated by automorphisms of the form $\exp(\ad(x_\a))$, $\a\in X\cup(-X)$, $x_\a\in \g_\a$.  The elements of $E_X$ stabilize $\g_X$ and fix $\h_X^\perp$ pointwise.

In what follows we will consider automorphisms $\xi$ of 
finite order of $\g$ so that for some subset $X$
of $\Pi$ of finite type the following conditions hold:
\begin{itemize}
\item[\rm (i)] $\xi = \nu \omega \mu$, where $\nu, \omega, \mu\in \Ag$
pairwise commute.
\item[\rm (ii)] $\nu\in \Aut(A)$ and $\nu(X) = X$
\item[\rm (iii)] $\omega$ is the Chevalley involution
\item[\rm (iv)] $\mu\in E_X$ and $\mu(\h_X) = \h_X$.
\end{itemize}

Suppose now that $\sg$ is of the second kind.  Thus, $\sg$ is
a semisimple automorphism of second the kind, and so, by 4.38
and 4.39 of \cite{KW}, $\sg$ is conjugate in $\Ag$ to an automorphism
of the form $\xi \rho$, where $\xi$ is an automorphism of $\g$ of 
finite order satisfying (i)--(iv) above (for some $X$), 
$\rho\in \Aut_\k(\g;\h)$ and $\xi\rho= \rho\xi$. 
(Actually more information about $\xi$ is given in \cite{KW},
but this is all we need here.)

By Lemma \ref{Lemma:LoopBasic}(b), we may assume that $\sg = \xi\rho$.  
Then $\rho$ has finite order.  Hence, 
increasing $m$ if necessary (which we can do by Lemma \ref{Lemma:Loopm}), 
we can assume that $\xi^m = \rho^m = \id$.
Then, by Theorem \ref{Thm:ErasingKM}, we have $\Lp(\g,\sg) \isom_R \Lp(\g,\xi)$.
But $p(\sg) = p(\xi)$ and so we can assume that $\sg = \xi$.

Now $p(\xi) = \tau$, where 
$$\tau = \nu \omega.$$
Thus, it remains to show that
$L(\g,\tau\mu) \isom_R \Lp(\g,\tau)$.
To prove this we perform a ``Cartan switch''. That is, we
choose a subalgebra $\jfrak$ of $\g$ so that
\begin{itemize}
\item[\rm (a)] $\mu\in \Aut_\k(\g;\jfrak)$
\item[\rm (b)] $\tau(\jfrak) = \jfrak$
\item[\rm (c)] $\jfrak$ is $\Ag$-conjugate to $\h$.
\end{itemize}
Once we have chosen such a subalgebra $\jfrak$,
we will have $L(\g,\tau\mu) \isom_R \Lp(\g,\tau)$ by
Lemma \ref{Lemma:LoopBasic}(b) and Theorem \ref{Thm:ErasingKM},
thereby completing this proof.

To choose $\jfrak$, we first let
$$\f = \g_X^\mu = \{x\in \g_X : \mu(x) = x\}.$$
Then, since $\mu|_{\g_X}$ is a semisimple automorphism of $\g_X$, $\f$
is a reductive subalgebra of $\g_X$ \cite[Corollaire \`a Proposition 12, \S 1, ${}^\circ 4$]{B}.

Next, since $\nu$ and $\omega$  stabilize $\g_X$ and commute with $\mu$, 
both $\nu$ and $\omega$ stabilize~$\f$.  Hence, 
$\tau$ stabilizes $\f$.  Therefore, $\tau|_\f$ is a semisimple automorphism of a reductive Lie algebra, and so there exists a Cartan subalgebra $\jfrak_X$ of $\f$ that is stabilized
by $\tau$ (see \cite[Theorem 4.5]{BoM} or \cite[Theorem 9]{P1}).  Let
$$\jfrak = \jfrak_X \oplus \h_X^\perp.$$

We conclude by checking the conditions (a), (b) and (c) above.  First $\mu$ fixes $\jfrak_X$
and $\h_X^\perp$ pointwise, and so (a) holds.  Also, $\tau$ stabilizes $\h_X$
and $\h$ and hence also $\h_X^\perp$ (since $\tau$ preserves the form $\fm$ restricted to $\h$).
Thus, (b) holds.  For (c), note that $\mu|_{\g_X}$ is an elementary automorphism of $\g_x$;
that is $\mu|_{\g_X}$ is the product of automorphisms of the form $\exp(\ad(x))$, where $x\in \g_X$
and $\ad_{\g_X}(x)$ is nilpotent.  Thus, by \cite[Theorem 11]{P1}, any Cartan subalgebra
of $\g_X^\mu$ is a Cartan subalgebra of $\g_X$.  In particular $\jfrak_X$ is a Cartan subalgebra of
$\g_X$.  Hence, by the conjugacy theorem for Cartan subalgebras of $\g_X$ \cite[p.28 and 48]{Sel},
there exists $\varphi\in E_X$ so that $\varphi(\jfrak_X) = \h_X$.  But $\varphi$
fixes $\h_X^\perp$ pointwise and so $\varphi(\jfrak) = \h$.
This proves (c) and the proof of the theorem is complete.
\endproof

We can now prove the 
second main theorem of this paper.

\begin{thm} 
\label{Thm:Main}
Let $\k$ be an algebraically closed field of characteristic 0,
and let  $\g = \g(A)$ be the Kac-Moody Lie algebra over $k$ 
determined by an indecomposable symmetrizable GCM $A$. Suppose that
$\sg_1$ and $\sg_2$ are automorphisms of $\g$ of period $m$. Then
$$\Lp(\g,\sg_1)\isom_\k \Lp(\g,\sg_2) \iff
p(\sg_1)\sim p(\sg_2).$$
\end{thm}

\beginproof
The implication ``$\Longrightarrow$'' is proved in Theorem \ref{Thm:Necessary}.
Conversely, suppose that $p(\sg_1)\sim p(\sg_2)$. 
Then,
$$
\renewcommand{\arraystretch}{1.25}
\begin{array}{llll}
\Lp(\g,\sg_1) 
&\isom_\k& \Lp(\g,p(\sg_1))&\mbox{by Theorem \ref{Thm:Project}}\\
&\isom_\k& \Lp(\g,p(\sg_2))&\mbox{by Lemma \ref{Lemma:LoopBasic}}\\
&\isom_\k& \Lp(\g,\sg_2)&\mbox{by Theorem \ref{Thm:Project}.}\quad\mbox{\Halmos}
\end{array}
$$
\bigskip

As a corollary of this theorem, we can obtain the corresponding theorem for loop
algebras of $\g'$.  Recall from Proposition \ref{Prop:Aut} that we have
$$\Aut_\k(\g')= \Autz_\k(\g') \rtimes \Out(A).$$
Let $p' = p'_\k : \Ag \to \Out(A)$
be the projection onto the second factor relative to this decomposition. 
Clearly we have
\beq
\label{ppp}
p'(\tau|_{\g'}) = p(\tau)
\eeq
for $\tau\in\Ag$.
(Recall we are identifying $\Out(A)$ as a subgroup of $\Ag$ and of $\Aut_\k(\g')$.)

\begin{cor} 
\label{Cor:Main}
Let $\k$ and $\g$ be as in Theorem \ref{Thm:Main}.
Suppose that
$\sg_1$ and $\sg_2$ are automorphisms of $\g'$ of period $m$. Then
$$\Lp(\g',\sg_1)\isom_\k \Lp(\g',\sg_2) \iff
p'(\sg_1)\sim p'(\sg_2).$$
\end{cor}

\beginproof
By Proposition \ref{Prop:Aut}, we can choose $\tau_1,\tau_2\in \Ag$ so that
$\tau_i|_{\g'} = \sg_i$ for $i=1,2$.  Then
$$
\renewcommand{\arraystretch}{1.25}
\begin{array}{llll}
\Lp(\g',\sg_1)\isom_\k \Lp(\g',\sg_2)
&\iff \Lp(\g,\tau_1)\isom_\k \Lp(\g,\tau_2)
&\mbox{by Theorem \ref{Thm:Equivalent}}\\
&\iff p(\tau_1) \sim p(\tau_2)
&\mbox{by Theorem \ref{Thm:Main}}\\
&\iff p'(\sg_1) \sim p'(\sg_2)
&\mbox{by (\ref{ppp})}\quad\mbox{\Halmos}
\end{array}
$$
\bigskip

\begin{rem} (a)  If $\g$ is finite dimensional, then
Theorem \ref{Thm:Main} follows from Proposition 8.5 of \cite{K1} 
and the Peterson-Kac conjugacy theorem for Cartan subalgebras of 
$\Lpg$.  

(b) If $\g$ is finite dimensional, Theorem \ref{Thm:Main} 
has also been proved using purely cohomological methods (along with
Theorem \ref{Thm:Compare} to relate $\k$-isomorphism and
$R$-isomorphism) in \cite[Proposition 10]{P1}.  It would be interesting to have
a purely cohomological proof in the general  case.
\end{rem}

\big
\noindent 
Department of Mathematical and Statistical Sciences, \newline
University of Alberta, \newline
Edmonton, Alberta, Canada T6G 2G1.  \newline
E-mail: ballison@math.ualberta.ca

\smallskip\noindent 
Department of Mathematics and Statistics,\newline
University of Saskatchewan, \newline
Saskatoon, Saskatchewan, Canada S7N 5E6\newline
E-mail: berman@math.usask.ca

\smallskip\noindent
Department of Mathematical and Statistical Sciences,\newline
University of Alberta,\newline
Edmonton, Alberta, Canada T6G 2G1\newline
E-mail: a.pianzola@ualberta.ca
\end{document}